\documentclass[a4paper,english,10pt]{article}
 \usepackage[applemac]{inputenc}
\usepackage{wrapfig}
\usepackage[english]{babel}
\usepackage{ulem}
\usepackage{epsf}
\usepackage[dvips]{epsfig}
\usepackage{subfig}
 \usepackage{graphicx}
\usepackage{hyperref}
\usepackage{amssymb}
\usepackage{amsmath}
\usepackage[nice]{nicefrac}
\usepackage{amsfonts}
\usepackage{dsfont}
\usepackage[T1]{fontenc}
\usepackage{wasysym}
\usepackage{amssymb}
\usepackage{stmaryrd}
\usepackage{aeguill}
\usepackage{mathrsfs}
\usepackage{csvsimple}
\usepackage[table]{xcolor}
\DeclareMathOperator*{\argmin}{arg\,min} 
\usepackage{tikz,tkz-tab}

\newcommand{\finpreuvem}{\quad\square}

\renewcommand{\epsilon}{\varepsilon}
\def\btab{\begin{eqnarray*}}
\def\etab{\end{eqnarray*}}
\def\beq{\begin{equation}}
\def\eeq{\end{equation}}

 \newcommand{ \un }{\mathds{1}}
 \newcommand{ \p }{\mathbb{P} }

 \newcommand{ \E }{\mathbb{E}}

 \newcommand{ \R }{ \mathbb{R} }
 \newcommand{ \Z }{ \mathbb{Z} }
 \newcommand{\N}{ \mathbb{N} }

\newcommand{ \bp}{ \mathbf{p}}

\newtheorem{The}{{\bf Theorem}}[section]
 \newtheorem{Lem}[The]{Lemma}
 \newtheorem{Cor}[The]{\bf Corollary}
 \newtheorem{Pro}[The]{\bf Proposition}
\newtheorem{Rem}[The]{{\bf Remark}}

 \newenvironment{proof}{\vspace{3mm}\noindent \textbf{Proof.}  }{\ \hfill $
 \finpreuvem$}

\makeatletter

\@addtoreset{equation}{section} \makeatother

\title{Power of the Crowd
\author{ A. Batakis, P. Debs, D. Le Peutrec  \footnote{Laboratoire IDP - C.N.R.S. UMR 7013 - Institut Denis-Poisson, Universit\'e d'Orl\'eans-Universit\'e de Tours
(France). \newline \vspace{0.1cm}  $\quad$  MSC 2000  . \newline \vspace{0.5cm} \textit{Key words: Galton Watson, Dynamical System, Elections, Random Walk. } }
}}

\begin{document}

\maketitle

\begin{abstract}
Consider a Galton Watson tree of height $m$: each leaf has one of $k$ opinions or not. In other words, for
$i\in \{1,\dots,k\}$,
 $x$ at generation $m$ thinks  $i$ with probability $\bp_i$ and nothing with probability~$\bp_0$. Moreover the opinions are independently distributed for
each leaf.\\
Opinions spread along the tree based on a specific rule: the majority wins! In this paper,  we study the asymptotic behavior of the distribution of the opinion of the root when $m\to\infty$.
\end{abstract}

\section{Introduction}
\begin{figure}
\begin{center}
\includegraphics[width=7cm]{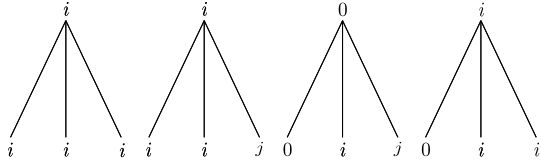}
\end{center}
\caption{The rules}
\label{fig1}
\end{figure}

\begin{figure}
\begin{center}
\includegraphics[width=7cm]{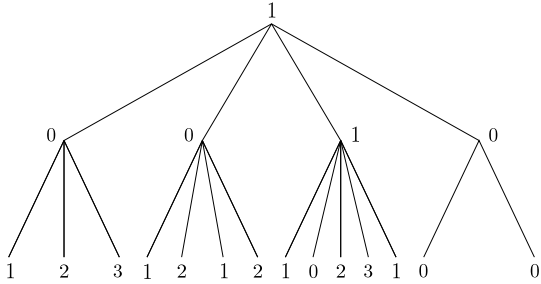}
\end{center}
\caption{An example}
\label{fig2}
\end{figure}

First let us recall the definition of a Galton Watson tree (GW) and give a few notations. Assume that $N$ is a $\mathbb N$-valued random variable following a distribution $q$: $\p(N=i)=q_i$ for $i\in\mathbb N$. In order to have a meaningful probabilistic setting, we assume that $q_0+q_1=0$ (B\"otcher case).\\ 
Let $\phi$ be the root of the tree and $ N_\phi$ an independent copy of  $ N$. Then, we draw $N_\phi$ children of~$\phi$: these individuals are the first generation. In the following we write $N$ for $N_\phi$ for typographical simplicity. At the $m$-th generation, for each individual $x$ we pick $ N_x$ an independent copy of $ N$ where $N_x$ is the number of children of $x$ and so on. The set $\mathbb T$, consisting of the root and its descendants, forms a GW of offspring distribution $q$.\\
We denote by $|x|$ the generation of $x$ and for $m\in\mathbb N$, $\mathbb T_m=\lbrace x\in \mathbb T,|x|\leq m\rbrace$ the GW cut at height~$m$ and the leaves of $\mathbb T_m$ are the elements of $\mathbb T_m\backslash \mathbb T_{m-1}$.  \\
Here we want to represent the propagation of an opinion in a population represented by a GW of height $m$.
More precisely consider  the set of probability vectors $\mathscr P_k$ defined by

\begin{equation}
\label{de.Pk}
\mathscr P_k:=\left\lbrace \bp=(\bp_0,\dots,\bp_k)\in(\mathbb R_+)^{k+1}: \sum_{i=0}^k\bp_i=1\ \ \text{and}\ \ \bp_{0}<1\right\rbrace
\subset\mathbb R^{k+1},
\end{equation}

and fix $\bp\in\mathscr P_k$. Each node of $\mathbb T_m$ has the opinion $\lbrace 1,\dots, k\rbrace$ according to the following rules: 
\begin{itemize}
\item Independently of the others, each leaf has an opinion according to $\bp$:
\[\p(\mbox{leaf thinks $i$})=\bp_i\ \ \ \text{and}\  \ \ \p(\mbox{leaf is undecided})=\bp_{0}.\]
\item The opinions spread to nodes at generation $m-1$ this way (see Figure \ref{fig1}):
\begin{enumerate}
\item[(R1)] the undecided children have no influence, except when the children are all undecided, in that case their ancestor has no opinion;    
\item[(R2)] if a relative majority of the children shares the same opinion, the ancestor thinks the same;
\item[(R3)] if several opinions are equally represented and the others are less, then the ancestor is undecided.
\end{enumerate}
\item We repeat this step for level $m-2$ and so on (see Figure  \ref{fig2}).
\end{itemize}

%

As claimed, we want to determine the asymptotic behavior of the distribution $\bp(m)$ of the state at the root of $\mathbb T_m$ when $m$ goes to infinity: 
\[\forall\, 1\leq i\leq k,\  \p(\mbox{root thinks }i)=\bp_i(m)\ \ \ \text{and}\ \ \ \p(\mbox{root is undecided})=\bp_{0}(m). \]
The children $\phi i$ of the root $\phi$ of $\mathbb T_{m}$, with $i\in\lbrace1,\dots,N\rbrace$,  are root nodes of $N$ independent GW of height $m-1$. Then the distribution $\bp(m)=\left(\bp_0(m),\bp_1(m),\dots,\bp_k(m)\right)$ of the state of $\phi$ is completely determined by the distribution  $\bp(m-1)$ of the independent states
of the $\phi i$,  $i\in\lbrace1,\dots,N\rbrace$. Let $H:\mathscr P_{k}\rightarrow \mathscr P_{k}$ be the function satisfying $\bp(m)=H(\bp(m-1))$, cf. \eqref{H}. An obvious reasoning by induction implies $H^m(\bp)=\bp(m)\in \mathscr P_{k} $. As a result, our problem consists in studying the orbits of  $H$ in $\mathscr P_{k}$. \\

The case of the binary tree is completely studied in \cite{Benjamini} and our paper can be seen as its natural generalization. Note that the relative majority is not the only possible extension of \cite{Benjamini}: in \cite{Haber}, the authors replaced (R2) and (R3) by the following rules 
\begin{enumerate}  
\item[(R2)'] if two children have different opinions, the ancestor is undecided;
\item[(R3)'] if all the children share the same opinion, the ancestor has it.
\end{enumerate}
We highlight the major differences between the results of \cite{Benjamini,Haber} and ours after the statement of our main results.\\

In what follows, we assume without loss of generality that $\bp_{1}>0$ and $\bp_1\geq\dots\geq\bp_k$.
It then holds $H_{1}(\bp)\geq \dots\geq H_{k}(\bp)$ and, if there exists $i\in\{1,\dots,k-1\}$
such that $\bp_{i+1}=\dots=\bp_{k}=0$, then $H_{i+1}(\bp)= \dots = H_{k}(\bp)=0$ (see remark \ref{monotonicity}).

It is hence sufficient 
to study the behavior of $H$ when acting on $\mathscr P_{i}$.
If there exists $i\in\lbrace1,\dots, k-1\rbrace$ such that  $\bp_1=\dots=\bp_i>\bp_{i+1}\geq\dots\geq\bp_{k}>0$, then ${i+1},\dots,k$ are called minor opinions and otherwise, i.e. if $\bp_1=\bp_2=\dots=\bp_k>0$, we say that we are in the uniform case.

In Section~\ref{se.2}, we prove that the major opinions do not vanish when $m\to\infty$, contrary to the minor opinions, and we state in Proposition \ref{nonunif_CV} a sufficient criterion to reduce the analysis to the uniform case. The biggest advantage of the uniform case is to study  the fixed points of a function defined on (a subset of) $\mathbb R$ instead of those of a function on $\mathscr P_{k}$. 
It naturally follows that if there is only one major opinion, regardless of the law of reproduction of $N$, this opinion spreads a.s. to the root asymptotically.\\ 
%
Although we have stated a very general problem, our main results below are available in  more restrictive cases: we only consider $n$-ary trees for $n\ge 2$ or GW trees supported in $2\mathbb N+1$ and two major opinions. This includes in particular a binary (``for-against") referendum, an election with two candidates. In the case of a $n$-ary tree for $n\ge2$, we obtain the following
\begin{The}\label{only}
For every $k\geq 2$ and $\bp\in \mathscr P_k$ such that $\bp_1=\bp_2>\bp_{3}\geq\dots\geq \bp_k$ and $\bp_{1}<\frac12$, $\bp(m)$ converges to $\left( \alpha_n,\frac{1-\alpha_n}{2},\frac{1- \alpha_n}{2},0_{k-2}\right)$ when $m\to\infty$, where $\alpha_n$ is the unique fixed point in~$(0,1)$ of the function
\[f_n:t\in[0,1]\longmapsto\sum_{k,0\le 2k\le n}\binom{n}{2k}\binom{2k}{k}\left(\frac{1-t}{2}\right)^{2k}t^{n-2k} \in[0,1].\]
Moreover, the above convergence remains true when $\bp_{1}=\bp_{2}=\frac12$ and $n$ is even,
whereas $\bp(m)=(0,\frac{1}{2},\frac{1}{2},0_{k-2})$ for every $m\geq 0$ when $\bp_1=\bp_2=\frac12$ and $n$ is odd.
\end{The}

The following result on the GW trees is ``just" a corollary as it needs to add a tricky argument to the proof of Theorem \ref{only} in the odd case.
\begin{Cor}\label{Corcool}
Taking a GW tree whose support is included in $2\mathbb N+1$ and such that $\E[N^2]<\infty$, the result of Theorem \ref{only} for odd $n$ remains true, replacing  $\alpha_n$ by $\alpha$, the unique fixed point in $(0,1)$ of 
\[f:t\in[0,1]\longmapsto\sum_{n\ge 1}q_{2n+1}f_{2n+1}(t).\]
\end{Cor}

As claimed, we make a brief list of the differences in \cite{Benjamini} and \cite{Haber}.\\
In \cite{Benjamini}, one can see that for $n=2$,  the result of Theorem \ref{only} is still true for any number of major opinions and the limit is explicit, in other words for all $i\ge 2$ and $\bp_1=\cdots=\bp_i>\bp_{i+1}\geq\dots\geq \bp_k$, $\bp(m)$ converges to  $\left(\frac{i-1}{2i-1},\frac{1}{2i-1},\dots, \frac{1}{2i-1},0_{k-i}\right)$.\\
With the rules explained in \cite{Haber}, the function studied in the uniform case with a $n$-ary tree and $i$ major opinions is the following:
\[g_n:t\in[0,\nicefrac{1}{i}]\longmapsto(1-(i-1)t)^n-(1-it)^n\]
and one can see this as the probability that the first opinion spreads. This function admits a unique fixed point $\bar x$ in $(0,\nicefrac{1}{i}]$ and the authors show that for $n\in\lbrace 3,4,5\rbrace$, $\bp(m)$ converges to $(1-i\bar x, \bar x,\dots,\bar x,0_{k-i})$. For $n\ge6$, stranger things happen: for instance for $n=6$, $\bar x$ is a repelling fixed point of $g_n$ and, if $i=2$, the authors show that there is a unique attracting orbit of prime period 2. Moreover, numerical simulations suggest the existence of a unique attracting orbit for every $n$ and $i$.\\

Let us come back to the organization of the paper: in  Section~\ref{se.3}, we give the proof of Theorem~\ref{only} and Corollary \ref{Corcool}. If the support of the GW tree is a subset of $2\mathbb N$, we have just succeeded to prove that if there is convergence, it does to $\left(\alpha,\frac{1-\alpha}{2},\frac{1-\alpha}{2},0_{k-2}\right)$, where  $\alpha$ is the unique fixed point in~$[0,1)$ of 
\[f:t\in[0,1]\longmapsto\sum_{n\ge 1}q_{2n}f_{2n}(t).\]
We have the attractivity but do not succeed to prove that we have convergence for all $\bp \in \mathscr P_k$. In this section, we also provide bounds for the values of the fixed points $\alpha_n$ of the functions $f_n$.\\
The general case seems to be unreachable for the moment, we just have proved the existence of a non-repulsive fixed point, not even its uniqueness. Nevertheless, we give an example where everything works, the geometric law, in Section~\ref{se.4}.\\
Finally, in Section~\ref{se.5} we make some remarks and give open questions and  Section~\ref{se.6} is an Appendix.

\section{Reduction to the Uniform case}
\label{se.2}

As stated in the introduction, the aim of the present paper is the study of a particular type of dynamical systems: more precisely, given the function $H:\mathscr P_{k}\rightarrow\mathscr P_{k}$ defined below (see \eqref{H}), we are interested in the behavior of the orbits $H^{\ell}(x)$
of the elements $x\in\mathscr P_{k}$ (see \eqref{de.Pk})  when $\ell$ goes to infinity.
In this section we give sufficient conditions to reduce the problem to a subfamily  of functions $H$  corresponding to the uniform case, namely the functions
$h_{k}$ defined below (see \eqref{de.hk}).



First, we need to specify the function $H$. Summing on the number of children of a ``typical'' node and on the number of children with a neutral opinion,  the probability that, for $i=1,\dots,k$, the $i$-th opinion spreads to their parent is equal to 
\begin{equation}\label{haf}
\sum_{n\ge 2} q_n\sum_{m_{0}=0}^{n-1}\binom{n}{m_{0}}\bp_{0}^{m_{0}}\sum_{S^i_{n-m_{0}}}\binom{n-m_{0}}{m_1,m_2,\dots, m_k}\prod_{j=1}^{k}\bp_j^{m_j}\,,
\end{equation}
where $S^i_{n}=\lbrace (m_1,\dots, m_k)\in \mathbb N^k,\forall j\ne i,  m_j<m_i\in\mathbb N,\sum_{j=1}^{k} m_j=n\rbrace$ and $\binom{n}{m_1,m_2,\dots, m_k}$ is the multinomial coefficient. 
Our problem then requires to study the fixed points
in $\mathscr P_{k}$
 of the function $H:\mathscr P_{k}\rightarrow \mathbb R^{k+1}$
defined by:
\begin{equation}\label{H}
H_i(\bp_0,\dots,\bp_k)=\left\lbrace\begin{array}{ccc}
\sum\limits_{n\ge2}q_n\sum\limits_{m_{0}=0}^{n-1}\binom{n}{m_{0}}\bp_{0}^{m_{0}}\sum\limits_{S^i_{n-m_{0}}}\binom{n-m_{0}}{m_1,m_2,\dots, m_k}\prod\limits_{j=1}^{k}\bp_j^{m_j} &\mbox{when $i\ne 0$,}\\
1-\sum\limits_{j=1}^kH_j(\bp_0,\dots,\bp_k).     &\mbox{when $i= 0$.}
\end{array} \right.
\end{equation}
\begin{Rem}\label{monotonicity}
Note that $\mathscr P_{k}$ is stable by $H$ and that, for $\bp\in \mathscr P_{k}$,
we can assume without loss of generality that  $\bp_1\geq\dots\geq\bp_k$ (which implies $\bp_{1}>0$
by definition of $\mathscr P_{k}$, see \eqref{de.Pk}).
\\
In this case,
it holds $H_{1}(\bp)\geq \dots\geq H_{k}(\bp)$
as well as, for every $i\in\{1,\dots,k\}$: $H_{i}(\bp)>0$ if, and only if, $\bp_{i}>0$.
\\
In particular, if there exists $i\in\{1,\dots,k-1\}$
such that $\bp_{i+1}=\dots=\bp_{k}=0$, then $H_{i+1}(\bp)= \dots = H_{k}(\bp)=0$ and  it is thus  sufficient 
to study $H:\mathscr P_{i}\subset \mathbb R^{i+1}\to \mathscr P_{i}$.
\end{Rem}

In what follows we denote, for $k\in \N$,   
\begin{equation}\label{Q}
\mathscr Q_{k}:=\{\bp\in \mathscr P_{k}\,,\ \bp_{1}\geq\dots\geq \bp_{k}\}
\end{equation}
and according to the previous remark we only need to consider the action of the function $H$ on~$\mathscr Q_{k}$.

In the uniform case, i.e. when $\bp_1=\dots=\bp_k\in(0,\frac1k]$,
one has simply
$H(1-k\bp_{1},\bp_1,\dots,\bp_1)=(1-kh_{k}(\bp_{1}),h_{k}(\bp_{1}),\dots,h_{k}(\bp_{1})) $, where~$h_{k}$ is the real function
defined on $[0,\frac1k]$ by
\begin{equation}
\label{de.hk}
h_k(x)=\sum_{n\ge2}q_n\sum_{m_{0}=0}^{n-1}\binom{n}{m_{0}}(1-kx)^{m_{0}}\sum_{S^1_{n-m_{0}}}\binom{n-m_{0}}{m_1,m_2,\dots, m_k}x^{n-m_{0}}\in\left[0,\frac1k\right]\,.
\end{equation}
The  study of the fixed points in $\mathscr Q_{k}$ of $H$ in the uniform case  is thus  reduced  to the study of the fixed points
 in $(0,\frac1k]$ of $h_{k}$.
Note also here that $0$ is a fixed point of $h_{k}$ which is repulsive, since
$$
h'_{k}(0)=\sum_{n\ge2}q_n\binom{n}{n-1}\sum_{S^1_{1}}\binom{1}{m_1,m_2,\dots, m_k}
=\sum_{n\ge2}n\,q_n=\E[N]\geq 2\,.
$$
Let us also recall that the generating function $G$ of $N$ is defined by
\begin{equation}
\label{de.gen}
\forall s\in[-1,1],\, G(s)=\E\left[s^N\right]=\sum_{n\geq0}s^n\p(N=n)
=\sum_{n\geq0}s^nq_{n}. 
\end{equation}
 On $(-1,1)$, $G$ is $\mathscr C^\infty$  and:
 \begin{equation}
\forall k\in\mathbb N,\, G^{(k)}(s)=\E[N(N-1)\dots(N-k+1)s^{N-k}]\,,
 \end{equation}
 implying that 
  \begin{equation}
\forall k\in\mathbb N,\, G^{(k)}(1^-)=\E[N(N-1)\dots(N-k+1)]\quad \mbox{ and }\quad G^{(k)}(0)=k!\,q_{k}.
 \end{equation}
 In particular, we have here $G(0)=\p(N=0)=0$, $G'(0)=\p(N=1)=0$,  and $G^\prime(1^-)=\E[N]\geq 2$.

\begin{Lem}\label{zeroplusbis}
Assume that $\bp\in\mathscr Q_{k}$  and that
 $G^{(2)}(1)$ is finite. Then, there exists $\eta>0$ such that 
\[\forall m\in \mathbb N,\ \  \bp_1(m)\geq \beta:=\min \lbrace \eta^{\mathfrak a}q_{\mathfrak a},\bp_1\rbrace\,,\]
where $\mathfrak{a}:=\inf\lbrace n\ge2, q_n\ne0\rbrace$. 
\end{Lem}

\begin{proof}
For $\bp\in\mathscr{Q}_k$, we have $\bp_1>0$, $\bp_{1}\geq \cdots\geq \bp_{k}$, and $\bp_0=1-\sum_{i=1}^k\bp_i$. We get
\begin{eqnarray}\label{PierreGenius}
\nonumber H_1(\bp)>\sum_{n\ge2}q_n\binom{n}{n-1}\bp_{0}^{n-1}\sum_{S^1_{1}}\binom{1}{m_1,m_2,\dots, m_k}\prod_{j=1}^{k}\bp_j^{m_j}
&=&\sum_{n\geq2} q_nn\bp_1\bp_{0}^{n-1} \\
 &=&\bp_1G^\prime\left(1-\sum_{i=1}^k\bp_i\right).
\end{eqnarray}

Since
 $G^{(2)}(1)=\E[N(N-1)]\in\mathbb R_{+}^{*}$, one can write 
\[G^\prime(1-t)=G^\prime(1)-tG^{(2)}(1)+\varepsilon(t),\]
where $\frac{\varepsilon(t)}{t}\underset{t\rightarrow0}{\longrightarrow}0$. 

As a result, there exists $0<\eta^\prime<\nicefrac{G^\prime(1)}{3G^{(2)}(1)}$ such that $|\varepsilon(t)|\le \frac{tG^{(2)}(1)}{2}$ when  $0\le t\le \eta^\prime$. Then, for $0\le x\le \eta:=\frac{\eta^\prime}{k}$ and $0\leq y \leq(k-1)x$,
\begin{align*}
xG^\prime(1-x-y)&\ge x\left(G^\prime(1)-(x+y)G^{(2)}(1)-(x+y)\frac{G^{(2)}(1)}{2}\right)=x\left(G^\prime(1)-3(x+y)\frac{G^{(2)}(1)}{2}\right)\\
&\ge x \left(G^\prime(1)-\frac{G^{\prime}(1)}{2}\right)=x\frac{G^{\prime}(1)}{2}\ge x, 
\end{align*}
where the last inequality follows from $G^\prime(1)=\E[N]\geq 2$. 

Thus, according to \eqref{PierreGenius}, $H_1(\bp)\ge \bp_1$
when $\bp_1\le \eta$.
In addition, when $\bp_1>\eta$ :
\begin{eqnarray*}
H_1(\bp)\ge \sum_{n\ge \mathfrak{a}}q_n\binom{n}{n-\mathfrak a}\bp_0^{n-\mathfrak a}\bp_1^{\mathfrak a}= \frac{\bp_1^{\mathfrak a}}{\mathfrak a!}G^{(\mathfrak a)}(\bp_0)\ge \frac{\bp_1^{\mathfrak a}}{\mathfrak a!}G^{(\mathfrak a)}(0)\ge  \eta^{\mathfrak a}q_{\mathfrak a}>0. 
\end{eqnarray*}
An obvious recurrence gives the claimed result.
\end{proof}

\begin{Rem}
\label{re.PG}
Applying the relation (\ref{PierreGenius}) to a fixed point $(\bp_0,\bp_{1},\dots,\bp_{k})\in\mathscr Q_{k}$
of $H$ we get
\[\bp_{1}=H_1(\bp)>\sum_{n\ge2}q_n\binom{n}{n-1}\bp_{0}^{n-1}\bp_1=\bp_1G^\prime(\bp_{0})\,,\] 
then, since $\bp_{1}>0$,
\[G^\prime(\bp_{0})<1.\]
\end{Rem}

The following elementary lemma ensures the validity of further results using the differentiability of $H$ on $\mathscr P_k$.
\begin{Lem}\label{derivability}
Suppose that $G^{\prime}(1)$ is finite. Let $\mathscr J_k:=\left\lbrace \bp=(\bp_0,\dots,\bp_k)\in(\mathbb R_+)^{k+1}: \sum_{i=0}^k\bp_i\le 1\right\rbrace$ and still denote by $H_i$  the functions defined by (\ref{H}) on  $\mathscr J_k$ for $1\le i\le k$. Then, these functions are of class $\mathscr C^1$ on $\mathscr J_k$.
\end{Lem}
\begin{proof}
Note first that by definition, for every $1\le i\le k$ and 
$\bp\in \mathscr J_k$, $H_i(\bp)$ writes:
\[H_i(\bp):= \sum\limits_{n\ge2}q_n H_{i,n}(\bp) :=\sum\limits_{n\ge2}q_n\sum_{m\in{\mathbb N}^{k+1},|m|=n}a_{n,m}\bp^m,\]
where $a_{n,m}\in \mathbb R^{+}$ and, for every $m=(m_{0},\dots,m_{k})\in {\mathbb N}^{k+1}$, 
$|m|:=\sum_{i=0}^{k}m_{i}$ and $\bp^{m}:=\prod_{i=0}^{k}\bp_{i}^{m_{i}}$.
More precisely, for every $n\geq 2$ and $m\in{\mathbb N}^{k+1}$ satisfying $|m|=n$,
\begin{equation*}
0\leq a_{n,m}\leq \binom{n}{m}
\quad \text{and hence} \quad
H_{i,n}(\bp) \leq \!\!\!\!
\sum_{m\in{\mathbb N}^{k+1},|m|=n}\binom{n}{m}\bp^m = 
(\bp_{0}+\cdots+\bp_{k})^{n} \leq\ 1.
\end{equation*}
It follows that, for every $1\le i\le k$,
$H_{i}$ is continuous on $\mathscr J_k$ and satisfies
$$
\forall \bp\in\mathscr J_k,\ \ \ 
H_i(\bp)=\sum\limits_{n\ge2}q_n H_{i,n}(\bp)
\leq 
\sum\limits_{n\ge2}q_n(\bp_{0}+\cdots+\bp_{k})^{n}=
 G(\bp_{0}+\cdots+\bp_{k}).
$$
Moreover,
 for every $\ell\in\{0,\dots,k\}$, $\bp\in \mathscr J_k$,  and $n\geq 2$:
 \begin{align*}
\label{eq.H-C1}
0\leq \frac{\partial H_{i,n}}{\partial \bp_{\ell}}(\bp)
\leq \frac{\partial (\bp_{0}+\cdots+\bp_{k})^{n}}{\partial \bp_{\ell}}
 = n(\bp_{0}+\cdots+\bp_{k})^{n-1}\leq n.
\end{align*}
This implies the claimed result, since $G^{\prime}(1)=\sum\limits_{n\geq 2}nq_{n}$ is finite.
%
\end{proof}

The following lemma ensures that the minor opinions can not spread to the root asymptotically: 
\begin{Lem}\label{zerobis}
Assume that
 $G^{(2)}(1)$ is finite.
In the (non uniform) case with
$i<k\in\mathbb N^{*}$ and
$\bp\in\mathscr Q_{k}$ such that
 $\bp_1=\dots= \bp_i>\bp_{i+1}\geq \dots\geq \bp_{k}\geq 0$, 
it holds $\bp_j(m)\underset{m\rightarrow\infty}{\rightarrow}0$ for every $j\in\{i+1,\dots,k\}$. 
\end{Lem}
\begin{proof}
Note that we just have to prove that $\lim\limits_{n\to\infty}\bp_{i+1}(n)=0$ when $\bp_{i+1}>0$. 
In this case,
writing $w_n=\frac{\bp_{i+1}(n)}{\bp_1(n)}>0$, we can easily see that  for every $ n\geq 0$, $w_{n+1}=w_n u_n$, where:
\begin{equation*}\label{dec}
u_n:=\frac{\displaystyle\sum_{z\ge2}q_z\sum_{m_{0}=0}^{z-1}\binom{z}{m_0}\bp_0^{m_0}(n)\sum_{S_{z-m_0}^{1}}\binom{z-m_0}{m_1,m_2,\dots, m_k}\bp_{i+1}^{m_1-1}(n)\bp_1^{m_{i+1}}(n)\prod_{j=2, j\neq i+1}^{k}\bp_j^{m_j}(n) }{\displaystyle\sum_{z\ge2}q_z\sum_{m_0=0}^{z-1}\binom{z}{m_0}\bp_0^{m_0}(n)\sum_{S_{z-m_0}^{1}}\binom{z-m_0}{m_1,m_2,\dots, m_k}\bp_1^{m_1-1}(n)\bp_{i+1}^{m_{i+1}}(n)\prod_{j=2,j\neq i+1}^{k}\bp_j^{m_j}(n)}.
\end{equation*}
For every $(m_1,\dots,m_k)\in S^{1}_{z-m_0}$, since $\bp_1> \bp_{i+1}>0$, we have that $\bp_{i+1}^{m_1-1}(n)\bp_1^{m_{i+1}}(n)< \bp_1^{m_1-1}(n)\bp_{i+1}^{m_{i+1}}(n)$ when
$m_{i+1}<m_{1}-1$,
 implying that $0<w_{n+1}<w_n$. Thus, $(w_n)$ is a positive decreasing sequence, and consequently converges
 to some $\ell\geq 0$. Since $w_0<1$, note that $\ell<1$.  \\
By compacity, there exists moreover a subsequence ${n_m}$ such that $\lim\limits_{m\to\infty}\bp_{j}(n_m)=a_j$ for every $j=0,\dots, k$. From 
Lemma~\ref{zeroplusbis}, we have $a_1>0$. Now, assume that $a_{i+1}>0$. Since $\ell<1$, we have $a_1>a_{i+1}>0$ and, using the definition of $u_n$:
\begin{eqnarray*}
\lim_{m\rightarrow\infty}u_{n_m}=\frac{\displaystyle\sum_{z\ge2}q_z\sum_{m_0=0}^{z-1}\binom{z}{m_0}a_{0}^{m_0}\sum_{S^1_{z-m_0}}\binom{z-m_0}{m_1,m_2,\dots, m_k}a_{i+1}^{m_1-1}a_1^{m_{i+1}}\prod_{j=2, j\neq i+1}^{k}a_j^{m_j} }{\displaystyle\sum_{z\ge2}q_z\sum_{m_0=0}^{z-1}\binom{z}{m_0}a_{0}^{m_0}\sum_{S^1_{z-m_0}}\binom{z-m_0}{m_1,m_2,\dots, m_k}a_1^{m_1-1}a_{i+1}^{m_{i+1}}\prod_{j=2,j\neq i+1}^{k}a_j^{m_j}}=\ell^\prime<1.
\end{eqnarray*}
This permits us to easily conclude:
\begin{eqnarray*}
\ell=\lim_{m\rightarrow\infty}w_{n_m+1}
=\lim_{m\rightarrow\infty}w_{n_m}u_{n_m}=\ell \ell^\prime
<\ell\,,
\end{eqnarray*}
which is a contradiction. Then $a_{i+1}=0$ and consequently $\bp_{i+1}(n)\underset{n\rightarrow\infty}{\rightarrow}0$. 
\end{proof}

\begin{Rem}
\label{re.expo0}
Actually,
for every $j\in\{i+1,\dots,k\}$,
 the convergence 
$\bp_j(m)\underset{m\rightarrow\infty}{\rightarrow}0$ 
is exponential, i.e. 
\begin{equation}
\label{eq.expo0}
\exists\, a \in(0,1)\,,\ \exists \,C >0\,,\ \forall \,m\in\N\,:\ \  0\leq \bp_{k}(m)\leq\dots\leq \bp_{i+1}(m)\leq C a^{m}\,.
\end{equation}
Note that to prove \eqref{eq.expo0}, it suffices to show that $\limsup\limits_{n\to\infty}u_{n}<1$, where
$(u_{n})_{n\geq 0}$ is the positive sequence introduced in the preceding proof. 
But since $\bp_{i+1}(m)\underset{m\rightarrow\infty}{\rightarrow}0$,
there exists a sequence $(\varepsilon_{n})_{n\geq0}$ converging to $0$ such that:
\begin{align*}
 u_n&=\frac{\displaystyle\sum_{z\ge2}q_z\binom{z}{z-1}\bp_0^{z-1}(n)\  +\  \varepsilon_{n} }{\displaystyle\sum_{z\ge2}q_z\sum_{m_0=0}^{z-1}\binom{z}{m_0}\bp_0^{m_0}(n)\sum_{S_{z-m_0}^{1}}\binom{z-m_0}{m_1,m_2,\dots, m_k}\bp_1^{m_1-1}(n)\bp_{i+1}^{m_{i+1}}(n)\prod_{j=2,j\neq i+1}^{k}\bp_j^{m_j}(n)}\\
&\leq 
\frac{\displaystyle\sum_{z\ge2}q_z\binom{z}{z-1}\bp_0^{z-1}(n)\  +\  \varepsilon_{n} }{\displaystyle
\sum_{z\ge2}q_z\binom{z}{z-1}\bp_0^{z-1}(n)+
\sum_{z\ge2}q_z\binom{z}{0}\bp_1^{z-1}(n)}
\\
&=\frac{G'(\bp_0(n))\  +\  \varepsilon_{n} }{
G'(\bp_0(n))+\bp_1^{-1}(n) G(\bp_1(n))}\ \leq\ 
\frac{G'(1-\beta) }{
G'(1-\beta)+G(\beta)}+  \frac{\varepsilon_{n}}{G(\beta)}
\ \underset{n\to\infty}{\longrightarrow}\ \frac{G'(1-\beta) }{
G'(1-\beta)+G(\beta)}<1\,,
\end{align*}
which implies $\limsup\limits_{n\to\infty}u_{n}<1$ and then \eqref{eq.expo0}.
\end{Rem}

In what follows, given a  real function $f$, we say that a fixed point $x$ of $f$ is linearly attracting for  $f$ when $f$ is differentiable at $x$ and $|f'(x)|<1$.


\begin{Pro}\label{dim1attrac}
Assume that
 $G^{\prime}(1)$ is finite.
Let $i\leq k\in\mathbb N^{*}$ 
and  assume that
 $\bar{x}_{i}\in(0,\frac1i]$ is a linearly attracting fixed point for the function $h_{i}$ defined in \eqref{de.hk}.
Then,  ${\bf{\bar x}}=(1-i\bar{x}_{i},\bar{x}_{i},\dots, \bar{x}_{i},0_{k-i})$ is an attracting fixed point for 
$$H\ :\ {\mathscr Q_{k,i}}
:=\{\bp\in\mathscr Q_{k}\,,\ \bp_{1}=\dots=\bp_{i}>\bp_{i+1}\geq\dots\geq \bp_{k}\}\to{\mathscr Q_{k,i}}\,.$$
\end{Pro}

\begin{proof} To prove that ${\bf{\bar x}}$ is attracting, note that it is sufficient to show that all the eigenvalues of the matrix $A:=\frac{\partial \tilde H}{\partial x}({\bf{y}})$ are in $(-1,1)$, where ${\bf{y}}=(\bar{x}_{i},0_{k-i})$ and $\tilde{H} = (\tilde{H}_1,\dots,\tilde{H}_{k-i+1}) = (H_1,H_{i+1},\dots, H_k)$ is a truncated version of $H$. For $\ell\in\{1,\dots,k-i+1\}$, $\tilde H_\ell$ is then defined by
\[
\sum_{n\ge 2} q_n\sum_{m_{0}=0}^{n-1}\binom{n}{m_{0}}\left(1-ix_1-\sum_{j=i+1}^kx_j\right)^{m_{0}}\sum_{S^{\tilde\ell}_{n-m_{0}}}\binom{n-m_{0}}{m_1,m_2,\dots, m_k}x_1^{\sum_{j=1}^{i}m_j}\prod_{j=i+1}^{k}x_j^{m_j}\,,
\]
where $\tilde \ell=1$ when $\ell=1$ and $\tilde\ell = \ell +i-1$ when $\ell\in \{2,\dots,k-i+1\}$.\\
Let us prove that the matrix $A$ is upper triangular, which
will immediately lead to the knowledge of its spectrum. For this purpose, let us compute 
$\frac{\partial \tilde H_\ell}{\partial x_{ r}}({\bf{y}})$ 
when $\tilde \ell\geq  r\in\{1,i+1,\dots,k\}$.\\
First, $\frac{\partial \tilde H_\ell}{\partial x_{r}}({{x_{1},x_{i+1},\dots,x_{k}}})$ 
equals, when $r=1$,
\begin{eqnarray*}
-\sum_{n\ge 2} q_n\sum_{m_{0}=0}^{n-1}\binom{n}{m_{0}}im_0\left(1-ix_1-\sum_{j=i+1}^kx_j\right)^{m_{0}-1}\sum_{S^{\tilde\ell}_{n-m_{0}}}\binom{n-m_{0}}{m_1,m_2,\dots, m_k}x_1^{\sum_{j=1}^{i}m_j}\prod_{j=i+1}^{k}x_j^{m_j}\\
+\sum_{n\ge 2} q_n\sum_{m_{0}=0}^{n-1}\binom{n}{m_{0}}\left(1-ix_1-\sum_{j=i+1}^kx_j\right)^{m_{0}}\sum_{S^{\tilde\ell}_{n-m_{0}}}\binom{n-m_{0}}{m_1,m_2,\dots, m_k}{\left(\sum_{j=1}^{i}m_j\right)}x_1^{\sum_{j=1}^{i}m_j-1}\prod_{j=i+1}^{k}x_j^{m_j}\,,
\end{eqnarray*}
and, when $r\in \{i+1,\dots,k\} $,
\begin{eqnarray*}
-\sum_{n\ge 2} q_n\sum_{m_{0}=0}^{n-1}\binom{n}{m_{0}}m_0\left(1-ix_1-\sum_{j=i+1}^kx_j\right)^{m_{0}-1}\sum_{S^{\tilde\ell}_{n-m_{0}}}\binom{n-m_{0}}{m_1,m_2,\dots, m_k}x_1^{\sum_{j=1}^{i}m_j}\prod_{j=i+1}^{k}x_j^{m_j}\\
+\sum_{n\ge 2} q_n\sum_{m_{0}=0}^{n-1}\binom{n}{m_{0}}\left(1-ix_1-\sum_{j=i+1}^kx_j\right)^{m_{0}}\sum_{S^{\tilde\ell}_{n-m_{0}}}\binom{n-m_{0}}{m_1,m_2,\dots, m_k}x_1^{\sum_{j=1}^{i}m_j}m_{ r} x_r^{m_r-1}\prod_{j=i+1, \neq r}^{k}x_j^{m_j}\,.
\end{eqnarray*}
Thus, by evaluating at ${\bf y}=(\bar{x}_{i},0_{k-i})$:
\begin{itemize}
\item When $\ell= r=1$,
\begin{eqnarray*}
\frac{\partial \tilde H_1}{\partial x_1}({{\bf y}})&=&-\sum_{n\ge 2} q_n\sum_{m_{0}=0}^{n-1}\binom{n}{m_{0}}im_0\left(1-i\bar x_i\right)^{m_{0}-1}\sum_{S^1_{n-m_{0}}}\binom{n-m_{0}}{m_1,m_2,\dots, m_k}{\bar x_i}^{\sum_{j=1}^{i}m_j}\prod_{j=i+1}^{k}0^{m_j}\\
&+&\sum_{n\ge 2} q_n\sum_{m_{0}=0}^{n-1}\binom{n}{m_{0}}\left(1-i\bar x_i\right)^{m_{0}}\sum_{S^1_{n-m_{0}}}\binom{n-m_{0}}{m_1,m_2,\dots, m_k}{\left(\sum_{j=1}^{i}m_j\right)}{\bar x_i}^{\sum_{j=1}^{i}m_j-1}\prod_{j=i+1}^{k}0^{m_j}\\
&=&-\sum_{n\ge 2} q_n\sum_{m_{0}=0}^{n-1}\binom{n}{m_{0}}im_0\left(1-i\bar x_i\right)^{m_{0}-1}\sum_{S^1_{n-m_{0}}}\binom{n-m_{0}}{m_1,\dots, m_i,0_{k-i}}{\bar x_i}^{n-m_{0}}\\
&+&\sum_{n\ge 2} q_n\sum_{m_{0}=0}^{n-1}\binom{n}{m_{0}}\left(1-i\bar x_i\right)^{m_{0}}\sum_{S^1_{n-m_{0}}}\binom{n-m_{0}}{m_1,\dots, m_i,0_{k-i}}(n-m_{0}){\bar x_i}^{n-m_{0}-1}\\
&=&h_{i}^\prime(\bar x_i)\,,
\end{eqnarray*}
where the function $h_{i}$ has been defined in \eqref{de.hk}.
\item When $\tilde\ell> r=1$,
\begin{eqnarray*}
\frac{\partial \tilde H_\ell}{\partial x_1}({{\bf y}})&=&-\sum_{n\ge 2} q_n\sum_{m_{0}=0}^{n-1}\binom{n}{m_{0}}im_0\left(1-i\bar x_i\right)^{m_{0}-1}\sum_{S^{\tilde\ell}_{n-m_{0}}}\binom{n-m_{0}}{m_1,m_2,\dots, m_k}{\bar x_i}^{\sum_{j=1}^{i}m_j}\prod_{j=i+1}^{k}0^{m_j}\\
&+&\sum_{n\ge 2} q_n\sum_{m_{0}=0}^{n-1}\binom{n}{m_{0}}\left(1-i\bar x_i\right)^{m_{0}}\sum_{S^{\tilde\ell}_{n-m_{0}}}\binom{n-m_{0}}{m_1,m_2,\dots, m_k}{\left(\sum_{j=1}^{i}m_j\right)}{\bar x_i}^{\sum_{j=1}^{i}m_j-1}\prod_{j=i+1}^{k}0^{m_j}\\
&=&0\,.
\end{eqnarray*}
\item Lastly, when $\tilde\ell\geq r>1$,
\begin{eqnarray*}
\frac{\partial \tilde H_\ell}{\partial x_r}({{\bf y}})&=&-\sum_{n\ge 2} q_n\sum_{m_{0}=0}^{n-1}\binom{n}{m_{0}}m_0\left(1-i\bar x_i\right)^{m_{0}-1}\sum_{S^{\tilde\ell}_{n-m_{0}}}\binom{n-m_{0}}{m_1,m_2,\dots, m_k}{\bar x_i}^{\sum_{j=1}^{i}m_j}\prod_{j=i+1}^{k}0^{m_j}\\
&+&\sum_{n\ge 2} q_n\sum_{m_{0}=0}^{n-1}\binom{n}{m_{0}}\left(1-i\bar x_i\right)^{m_{0}}\sum_{S^{\tilde\ell}_{n-m_{0}}}\binom{n-m_{0}}{m_1,m_2,\dots, m_k}{\bar x_i}^{\sum_{j=1}^{i}m_j}m_r 0^{m_r-1}\prod_{j=i+1, \neq r}^{k}0^{m_j}\\
&=&\sum_{n\ge 2} q_n\sum_{m_{0}=0}^{n-1}\binom{n}{m_{0}}\left(1-i\bar x_i\right)^{m_{0}}\sum_{S^{\tilde\ell}_{n-m_{0}}}\binom{n-m_{0}}{m_1,m_2,\dots, m_k}{\bar x_i}^{\sum_{j=1}^{i}m_j}m_r 0^{m_r-1}\prod_{j=i+1, \neq r}^{k}0^{m_j}
\end{eqnarray*}
which, when $\tilde\ell> r$, equals 0  and, when  $\tilde\ell= r$, equals
\begin{eqnarray*}
\sum_{n\ge 2} q_n\sum_{m_{0}=0}^{n-1}\binom{n}{m_{0}}\left(1-i\bar x_i\right)^{m_{0}}
\binom{n-m_{0}}{0_{\tilde\ell-1}, 1, 0_{k-\tilde\ell}}
=\sum_{n\ge 2} nq_n(1-i\bar x_i)^{n-1}
=G^\prime(1-i\bar x_i).
\end{eqnarray*}
\end{itemize}
As claimed, $A$ is  thus upper triangular and its spectrum is $\lbrace h_{i}^\prime(\bar x_i),G^\prime(1-i\bar x_i)\rbrace$. 
Moreover, $h'_{i}(\bar{x}_{i})$ belongs to $(-1,1)$ by assumption 
and, according to Remark~\ref{re.PG}, $G^\prime(1-i\bar x_i)$
also belongs to $(-1,1)$
(since $H({\bf{\bar x}})={\bf{\bar x}}\in\mathscr Q_{k}$). The statement of Proposition~\ref{dim1attrac} follows.
%
\end{proof}

The following proposition is an adaptation of Proposition 3.11 in \cite{Haber}.

\begin{Pro}\label{nonunif_CV}
Assume that
 $G^{(2)}(1)$ is finite.
Let $i\leq k\in\mathbb N^{*}$ 
and  assume that
 $\bar{x}_{i}\in[0,\frac1i]$ is a linearly attracting fixed point for the function $h_{i}$ defined in \eqref{de.hk}
 whose basin of attraction contains~$(0,\frac1i]$.
Then, $\bar{x}_{i}\in(0,\frac1i]$ and
${\bf{\bar x}}=(1-i\bar{x}_{i},\bar{x}_{i},\dots, \bar{x}_{i},0_{k-i})$ is a globally attracting fixed point for 
$$H\ :\ {\mathscr Q_{k,i}}
=\{\bp\in\mathscr Q_{k}\,,\ \bp_{1}=\dots=\bp_{i}>\bp_{i+1}\geq\dots\geq \bp_{k}\}\to{\mathscr Q_{k,i}}\,.$$ 
\end{Pro}


\begin{proof}
Note first that when
 $k\geq i$, $\bp_1=\dots=\bp_i>0$, and $\bp_j=0$ for 
 $j=i+1,\dots,k$, then $\bp(m)$  converges to $(1-i\bar{x}_{i},\bar{x}_{i},\dots, \bar{x}_{i}, 0_{k-i})$
 by hypothesis.
 According to Lemma~\ref{zeroplusbis}, it thus holds~$\bar{x}_{i}>0$.\\
We have now to extend this result  when $k> i$, $\bp_1=\dots=\bp_i>\bp_{i+1}\geq\dots\geq \bp_{k}\geq0$, and~$\bp_{i+1}>0$.
Let us again consider
the truncated version of $H$, 
 $\tilde{H} = (\tilde{H}_1,\dots,\tilde{H}_{k-i+1}) = (H_1,H_{i+1},\dots, H_k)$, where,
 for $\ell\in\{1,\dots,k+i-1\}$,  $\tilde H_\ell$ is defined by
\[
\sum_{n\ge 2} q_n\sum_{m_{0}=0}^{n-1}\binom{n}{m_{0}}\left(1-ix_1-\sum_{j=i+1}^kx_j\right)^{m_{0}}\sum_{S^{\tilde\ell}_{n-m_{0}}}\binom{n-m_{0}}{m_1,m_2,\dots, m_k}x_1^{\sum_{j=1}^{i}m_j}\prod_{j=i+1}^{k}x_j^{m_j}\,,
\]
where $\tilde \ell=1$ when $\ell=1$ and $\tilde\ell = \ell +i-1$ when $\ell\in \{2,\dots,k-i+1\}$.
Let us also define the set 
$$\tilde{\mathscr Q}_{k,i}:=\left\lbrace (x_1,\dots,x_{k-i+1})\in\mathbb R^{k-i+1},\frac1i\geq x_1>x_2\geq x_3\geq\dots \geq x_{k-i+1}\geq 0, ix_1+\sum_{j=2}^{k-i+1}x_j \leq1 \right\rbrace.$$ 
Let us show that, for every ${\bp} \in \tilde{\mathscr Q}_{k,i}$, $\tilde{H}^m({\bp})$ converges to ${\bf\bar{x}}= (\bar{x}_{i}, 0_{k-i})$,
which is equivalent to the convergence result stated
in Proposition~\ref{nonunif_CV}.
We fix $\bp=(\bp_1,\bp_{i+1},\bp_{i+2},\dots,\bp_{k})\in\tilde{\mathscr Q}_{k,i}$ and recall
from Lemma~\ref{zeroplusbis}
 that, for every $m\in\mathbb N$, $\big(\tilde{H}^m({\bp})\big)_{1} =\bp_{1}(m)\geq \beta:=\min \lbrace \eta^{\mathfrak a}q_{\mathfrak a},\bp_1\rbrace$.

As $\bar x_i$ is a linearly attracting fixed point of $h_i$,
for every $\varepsilon>0$ small enough, $B(\bar{x}_{i},\nicefrac{\varepsilon}{2})$ is $h_i$-invariant.
Now, while noting that $\tilde{H}^m (x,0_{k-i}) = (h_{i}^m(x),0_{k-i})$ for every $m\in\mathbb N$, we define $E_m:=\lbrace x\in\mathbb [0,1]: h_{i}^m(x)\in B(\bar{x}_{i},\nicefrac{\varepsilon}{2})\rbrace$
for some arbitrarily small $\varepsilon>0$. The sequence $(E_m)_{m\geq0}$ is an ascending chain of sets and from the convergence in the uniform case, 
\[[\beta,\nicefrac{1}{i}]\subset \bigcup_{m\geq 0} E_m. \]
As the inverse image of an open set of $\mathbb R$ by a continuous function, $E_m$ is an open set for all~$m\in\mathbb N$. Since $(E_m)_{m\geq0}$ is an increasing sequence of open sets covering the compact $[\beta,\nicefrac{1}{i}]$, there exists~$N\in\mathbb N$ such that  
\[[\beta,\nicefrac{1}{i}]\subset \bigcup_{m= 0}^N E_m = E_{N}\,, \]
implying that: $\forall x\in[\beta,\nicefrac{1}{i}], \, \tilde{H}^N(x,0_{k-i})\in B(\bar{x}_{i},\nicefrac{\varepsilon}{2}) \times \{0\}^{k-i} \subset B({\bf{\bar{x}}},\nicefrac{\varepsilon}{2})$. 

On the closed bounded set $\mathscr  G:=[\beta,\nicefrac{1}{i}]\times \mathbb R_+^{k-i}\cap \overline{\tilde{\mathscr Q}_{k,i}}$, $\tilde{H}^N$ is uniformly continuous and thus there exists $\delta>0$ such that 
\[\forall (x,y),(x^\prime,y^\prime)\in\mathscr G, \Vert (x,y)-(x^\prime,y^\prime)\Vert \leq \delta \Rightarrow \Vert \tilde{H}^{N}(x,y)-\tilde{H}^N(x^\prime,y^\prime) \Vert\leq \nicefrac{\varepsilon}{2}. \]  
According to Lemma~\ref{zerobis},  $\bp_{j}(m)=\big(H^m(\bp)\big)_{j}\rightarrow 0$ for every $j\in\{i+1,\dots,k\}$. Consequently, there exists $N_1\in\mathbb N$ such that, for every $m\in\mathbb N$:  $m\geq N_1$
implies
 $ \Vert (\bp_{i+1}(m),\dots, \bp_k(m))\Vert\leq \delta$ and then
\[  
\Vert \tilde{H}^{N}(\bp_1(m),\bp_{i+1}(m),\dots,\bp_{k}(m))-\tilde{H}^N(\bp_1(m),0_{k-i}) \Vert\leq \nicefrac{\varepsilon}{2}.\]
Thus, for every  $m\geq N_{1}$, the fact that $\bp_1(m)\in[\beta,\nicefrac{1}{i}]$ implies: 
\begin{align*}
 A:&=\Vert(\bp_1(m+N),\bp_{i+1}(m+N),\dots,\bp_{k}(m+N))-{\bf{\bar x}}\Vert\\
& =\Vert \tilde{H}^{N}(\bp_1(m),\bp_{i+1}(m),\dots,\bp_{k}(m))-{\bf{\bar x}}\Vert \\
&\leq \Vert \tilde{H}^{N}(\bp_1(m),\bp_{i+1}(m),\dots,\bp_{k}(m))-\tilde{H}^{N}(\bp_1(m),0_{k-i})\Vert
+\Vert \tilde{H}^{N}(\bp_1(m),0_{k-i})-{\bf{\bar x}}\Vert  \\
&\leq
 \nicefrac{\varepsilon}{2}+\nicefrac{\varepsilon}{2}=\varepsilon,
\end{align*}
which concludes the proof of Proposition~\ref{nonunif_CV} since  $\varepsilon>0$ is arbitrarily small.
\end{proof}


\begin{Rem}\label{light}
\begin{enumerate}
\item In the statement of Proposition \ref{nonunif_CV}, we can actually lighten the hypothesis $|h_i'(\bar{x}_i)|<1$
using the fact that $h_i$ is $\mathscr C^1$ on $[0,\frac1i]$ (which follows from Lemma~\ref{derivability}) and a monotonicity argument. However,
this would require unnecessary extra work as  this hypothesis  is satisfied
in all our examples.
\item It is not difficult to prove that if we have just one major opinion, it spreads almost surely to the root. Indeed, in the ``uniform'' case with only one opinion, according to the rules, the probability that in a GW of height $m$ the unique opinion does not spread to the root is equal to:  
\[\bp_0(m)=\sum_{n\ge 2}q_n\bp_0^n(m-1)=G(\bp_0(m-1))=G^{m}(\bp_0).\]
Since $N\ge 2$ a.s., $G$ is strictly convex  on $[0,1]$ with
$0$ and $1$ as
sole fixed points. It follows that :
$$
\forall \,\bp_0\in[0,1)\,,\ \ \lim_{m\rightarrow\infty}G^m(\bp_{0})=0\ \ \ \text{and thus}
\ \ \ 
\forall \,\bp_1\in(0,1]\,,\ \ \lim_{m\rightarrow\infty}\bp_1(m)=1\,.
$$
Proposition~\ref{nonunif_CV} then ensures the convergence in the non uniform case with one major opinion.
\end{enumerate}
\end{Rem}

\noindent
{\bf Conclusion :} 
From the above results, one deduces that:
\begin{itemize}
\item For any $\bp\in \mathscr Q_{k}$, defining
 $i:=\max\big\{\ell\in\{1,\dots,k\}\,, \bp_{\ell}=\bp_{1}\big\}\in\{1,\dots,k\}$,
we have $\bp_{j}(m)\underset{m\rightarrow\infty}{\rightarrow}0$ for every $j\in\{i+1,\dots,k\}$.\\
The  accumulation points of the sequence $(H^{\ell}(\bp))_{\ell\geq 0}$
have thus the form $(1-i\bar{x}_{i},\bar{x}_{i},\dots, \bar{x}_{i},0_{k-i})$
where, according to Lemma~\ref{zeroplusbis}, $\bar{x}_{i}\in(0,\frac1i]$.\\
In particular, the fixed points (resp. the $m$-cycles) of $H$ in $\mathscr Q_{k}$
are the $(1-i\bar{x}_{i},\bar{x}_{i},\dots, \bar{x}_{i},0_{k-i})$, where
$i\in\{1,\dots,k\}$ and $\bar{x}_{i}$ is a fixed point  (resp. a $m$-cycle) of $h_{i}$ in $(0,\frac1i]$.
\item Recall that, for $i\leq k\in \mathbb N^{*}$, 
$\mathscr Q_{k,i}=\{\bp\in \mathscr Q_{k}\,,\ \bp_{1}=\dots= \bp_{i}>\bp_{i+1}\}$
and $H(\mathscr Q_{k,i})\subset \mathscr Q_{k,i}$.\\
Proposition~\ref{dim1attrac} implies that the fixed point
$(1-i\bar{x}_{i},\bar{x}_{i},\dots, \bar{x}_{i},0_{k-i})$
of $H:\mathscr Q_{k,i}\to \mathscr Q_{k,i}$ is  attracting when
$\bar{x}_{i} \in (0,\frac1i]$ is linearly attracting  for $h_{i}$. Conversely, if
 $(1-i\bar{x}_{i},\bar{x}_{i},\dots, \bar{x}_{i},0_{k-i})$ is  attracting
for $H:\mathscr Q_{k,i}\to \mathscr Q_{k,i}$, then $\bar{x}_{i}\in(0,\frac1i]$
is obviously attracting for  $h_{i}$.\\
Finally, according to  Proposition~\ref{nonunif_CV}, if the basin of attraction of 
a fixed point $\bar{x}_{i} $ of $h_{i}$ in~$(0,\frac1i]$ is  $(0,\frac1i]$,
then the basin of attraction of $(1-i\bar{x}_{i},\bar{x}_{i},\dots, \bar{x}_{i},0_{k-i})$ for $H:\mathscr Q_{k}\to \mathscr Q_{k}$  is~$\mathscr Q_{k,i}$,
and the converse is clearly true.
\end{itemize}

\section{The 2 major opinions case or the second run of an election}
\label{se.3}

In this section, we consider only two major opinions.
Moreover, contrary to the previous section, we study the probability that the ``neutral'' opinion spreads, i.e. that in a group of $n$ individuals, no opinion has a majority. 
With in mind the results of the preceding section, we focus on the uniform case:
if $t\in[0,1]$ is the probability that a given individual gives a white vote, the probability
of each opinion is $\frac{1-t}{2}$ and the probability
 of the group to come up undecided is then given by 
\[H_{0}\left(t,\frac{1-t}{2},\frac{1-t}{2}\right)\ =\ 1-2h_{2}\left(\frac{1-t}{2}\right)\ :=\ f_n(t)=\sum_{k,0\le 2k\le n}\binom{n}{2k}\binom{2k}{k}\left(\frac{1-t}{2}\right)^{2k}t^{n-2k}.\]
We will thus study the fixed points of $f_{n}$ in $[0,1)$, or equivalently the fixed point of $h_{2}$ in $(0,\frac12]$.\\

We start by a crucial remark providing an integral formula for the functions $f_n$.
\begin{Lem}\label{LucDebs}
For all $0\le t\le1$:
\begin{equation}\label{Luc}
f_n(t)=\frac{1}{\pi}\int_{0}^\pi\left((1-t)\cos x+t\right)^n\mathrm{d}x.
\end{equation}
\end{Lem}

\begin{proof}
Recall the Wallis integral for all $k\ge0$:
\begin{equation}\label{Luc1}
\int_0^{\frac{\pi}{2}}\cos ^{2k} x\mathrm{d}x=\frac{\pi}{2}\frac{(2k)!}{(k!2^k)^2 }=\frac{\pi}{2^{2k+1}}\binom{2k}{k}\Leftrightarrow\frac{1}{2\pi}\int_{-\pi}^{\pi}\cos ^{2k}x\mathrm{d}x=\frac{1}{2^{2k}}\binom{2k}{k}
\end{equation}
and note that, using the substitution $u=\frac{\pi}{2}-x$ :
\begin{equation}\label{Luc2}
\int_{0}^{\pi}\cos ^{2k+1}x\mathrm{d}x=\int_{-\frac{\pi}{2}}^{\frac{\pi}{2}}\sin^{2k+1}u\mathrm{d}u=0.
\end{equation}
Then, using \eqref{Luc1} and \eqref{Luc2}, we can write:
\begin{align*}
f_n(t)&=\sum_{k,0\le 2k\le n}\binom{n}{2k}\left({1-t}\right)^{2k}t^{n-2k}\frac{1}{2^{2k}}\binom{2k}{k}\\
&=\sum_{k,0\le 2k\le n}\binom{n}{2k}\left({1-t}\right)^{2k}t^{n-2k}\frac{1}{2\pi}\int_{-\pi}^{\pi}\cos ^{2k}x\mathrm{d}x\\
&=\sum_{k,0\le k\le n}\binom{n}{k}\left({1-t}\right)^{k}t^{n-k}\frac{1}{2\pi}\int_{-\pi}^{\pi}\cos ^{k}x\mathrm{d}x\\
&=\frac{1}{2\pi}\int_{-\pi}^{\pi}\sum_{k,0\le k\le n}\binom{n}{k}\left({1-t}\right)^{k}t^{n-k}\cos ^{k}x\mathrm{d}x=\frac{1}{2\pi}\int_{-\pi}^\pi((1-t)\cos x+t)^n\mathrm{d}x, 
\end{align*}
and we easily conclude using the parity.
\end{proof}

There is a more elegant way to prove Lemma~\ref{LucDebs} by using Fourier series:
let us consider the random walk on $L^2({\mathcal S^1})$, the space of square integrable functions on the circle, defined on its usual $(e_{k}=e^{i k x})_{k\in\Z}$ basis by
$$Z_0=1\ \ \mbox{ and }\ \ \p(Z_{n+1}=e_k|Z_n=e_{\ell})=\frac{1-p}{2}\un_{k=\ell-1}+\frac{1-p}{2}\un_{k=\ell+1}+p\un_{k=\ell},$$
where $0\le p\le 1$.\\
Using that $ \p\left(\nicefrac{Z_{n+1}}{Z_n}=e^{ix}\right)=\p\left(\nicefrac{Z_{n+1}}{Z_n}=e^{-ix}\right)=\frac{1-p}{2}$, that $\p\left(\nicefrac{Z_{n+1}}{Z_n}=1\right)=p$, and the independence of the random variables $\nicefrac{Z_{n+1}}{Z_n}$, we get $$f_n(p)	=\p(Z_n=1).$$
The (infinite) matrix associated to the walk is 
$$ A=\left(\begin{array}{cccccccc}
\vdots& \ddots&\ddots& \ddots&\ddots& \ddots&\cdots&\vdots\\
0 &\cdots & p &  \frac{1-p}{2} & 0 &0  &\cdots & 0\\
0 &\cdots & \frac{1-p}{2} & p &  \frac{1-p}{2}&0& \cdots & 0\\
0 &\cdots & 0 & \frac{1-p}{2} & p &  \frac{1-p}{2}& \cdots & 0\\
0 &\cdots & 0 & 0& \frac{1-p}{2} & p & \cdots & 0\\
\vdots& \cdots&\ddots& \ddots&\ddots& \ddots&\ddots&\vdots
\end{array}\right),$$
so that $A$ applied to $e_{\ell}$ equals $\frac{1-p}{2}e_{\ell-1}+\frac{1-p}{2}e_{\ell+1}+pe_{\ell }$. Let $L$ be the associated linear operator on~$L^2({\mathcal S^1})$.
A straightforward easy computation shows that $ L(e^{i\ell x})=\left((1-p)\cos{x}+p\right)e^{i\ell x}$, which implies that $L$ is a scalar operator : 
$$L:L^2({\mathcal S^1})\ni h\mapsto \left((1-p)\cos{x}+p\right) h\in L^2({\mathcal S^1}) $$ and therefore the iterated operator $L^n$ is given by
$$L^n(h)= \left((1-p)\cos{x}+p\right)^n h.$$
On the other hand,
$$\p(Z_n=1)=<A^ne_0,e_0>=<L^n{\mathbf 1},{\mathbf 1}>=\frac1{2\pi} \int_0^{2\pi}\left((1-p)\cos{x}+p\right)^n{\rm d}x.$$

\begin{Rem}\label{deriv1}
As a polynomial of degree $n$, $f_n\in \mathscr C^{\infty}([0,1])$ and, for every $k\in\{0,\dots, n\}$:
\begin{equation}\label{deriv}
f_n^{(k)}(t)=\frac{1}{\pi}\frac{n!}{(n-k)!}\int_0^\pi(1-\cos x)^k(t(1-\cos x)+\cos x)^{n-k}\mathrm{d}x.
\end{equation}
\end{Rem}

\begin{figure}[h!]
\begin{center}
\includegraphics[width=6cm]{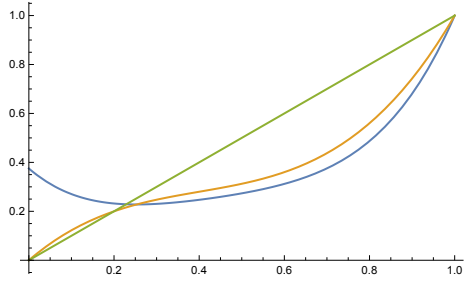}
\end{center}
\caption{The graphs of $f_n$ for $n=4$ (blue), $n=3$ (orange)}
\label{infcut}
\end{figure}

\begin{Lem}
\label{le.fix}
For all $n\ge 2$, the function $f_n$ admits in $(0,1)$ a unique fixed point $\alpha_n$, and
$\alpha_n<\nicefrac{1}{2}$.
\end{Lem}
\begin{proof}
To prove the unicity, we have to distinguish two cases according to the parity of $n$.
\begin{itemize}
\item \underline{Odd case (see the orange graph of Figure \ref{infcut}):}\\
Using Lemma \ref{LucDebs}, Remark \ref{deriv1} and \eqref{Luc1}
, we have:
\begin{align}
f_n(0)&=\frac{1}{\pi}\int_{0}^\pi\cos^n x \mathrm{d}x=0\ \ \text{and}\ \  f_n(1)=\frac{1}{\pi}\int_0^{\pi}\mathrm{d}x=1,\label{obvious}\\
f_n^{\prime}(0)&=\frac{n}{\pi}\left(\int_{0}^\pi\cos^{n-1} x\mathrm{d}x-\int_{0}^\pi\cos^{n} x\mathrm{d}x\right)=\frac{n}{2^{n-1}}\binom{n-1}{\frac{n-1}{2}}>1,\label{obvious2}\\
f_n^{\prime}(1)&=\frac{n}{\pi}\int_0^{\pi}(1-\cos x)\mathrm{d}x=n>1\label{obvious3},
\end{align}
implying that $f_n$ has at least one fixed point in $(0,1)$. 
The inequality in \eqref{obvious2} follows from (see \eqref{Nass1})
\begin{equation}
\label{de.xi-n}
\forall n \geq 1,\ \ \ 
{   \frac{2}{\sqrt{2\pi(2n+1)}}}< \xi_{2n}:=\frac{1}{2^{2n}}\binom{2n}{n}
\end{equation}
Note also that the formulas \eqref{obvious} are direct with the spreading rules.\\
Moreover, Remark  \ref{deriv1} with $k=3$ gives (note that $n\geq 3$ since $n$ is odd):
\[\forall t\in[0,1],\ \ f_n^{(3)}(t)=\frac{1}{\pi}\frac{n!}{(n-3)!}\int_0^\pi(1-\cos x)^3(t(1-\cos x)+\cos x)^{n-3}\mathrm{d}x> 0,\] 
implying that $f_n$ has at most three fixed points in $[0,1]$. As a result, $f_n$ has a unique fixed point $\alpha_n$ in $(0,1)$.
\item \underline{Even case (see the blue graph of Figure \ref{infcut}):}\\
Since  $f_n(0)=\frac{1}{2^{n}}\binom{n}{\frac{n}{2}}>0$, $f_n(1)=1$ and $f_n^\prime(1)=n>1$, we deduce that $f_n$ has at least one fixed point in $(0,1)$.
Using again Remark \ref{deriv1},   $f_n$ is strictly convex in $[0,1]$  and has thus  at most two fixed points in this interval. As a result,  $f_n$ has a unique fixed point $\alpha_n$ in $(0,1)$.
\end{itemize}
Thanks to the unicity of $\alpha_n$, we just have to show $f_{n}(\nicefrac{1}{2})<\nicefrac{1}{2}$ to obtain $\alpha_n<\nicefrac{1}{2}$.\\ 
According to the formula \eqref{luc1} of Lemma \ref{wiard3}:
\begin{equation}
\label{de.xi-n'}
f_n\left(\frac{1}{2}\right)=\frac{1}{2^n}\sum_{k,0\le 2k\le n}2^{-2k}\binom{n}{2k}\binom{2k}{k}=\frac{1}{2^{2n}}\binom{2n}{n}=\xi_{2n}.
\end{equation}
As $\xi_{2}=\nicefrac{1}{2}$ and $\left(\xi_{2n}\right)_{n\ge0}$ is a strictly decreasing sequence according to Lemma~\ref{XI},
$\xi_{2n}<\nicefrac{1}{2}$ for all $n\ge2$. 
\end{proof}

\begin{Rem}
\label{re.GW1}
If we look at the GW case, we have to study the fixed points in $[0,1)$ of:
\[f:t\in[0,1]\mapsto\sum_{n\ge 2}q_nf_n(t).\]
With similar arguments as those of the previous proof, it is not difficult to prove the existence of a fixed point $\alpha\in (0,\nicefrac{1}{2})$, since
\begin{align*}
f(1)=1,\,  f^\prime(1)&=\E[N]>1,\ \ \text{and}\ \   f(\nicefrac{1}{2})<\nicefrac{1}{2}.
\end{align*}
Indeed, if $f(0)=\sum_{n\ge 1}q_{2n}f_{2n}(0)>0$, we have our result and, otherwise,
$q_{2n}=0$ for every $n\geq 1$, so
 $f^\prime(0)=\sum_{n\ge1}q_{2n+1}f^\prime_{2n+1}(0)>1$ and we can easily conclude.
 
Moreover, if the support of $N$ is a subset of $2\mathbb N$ or of $2\mathbb N+1$, we have the unicity of $\alpha$
by the arguments used in the previous proof.
\end{Rem}

\subsection{The odd case}

\subsubsection{Basin of attraction of the fixed point {\boldmath $\alpha_{n}$}}

\begin{figure}[h!]
\begin{center}
\includegraphics[width=7cm]{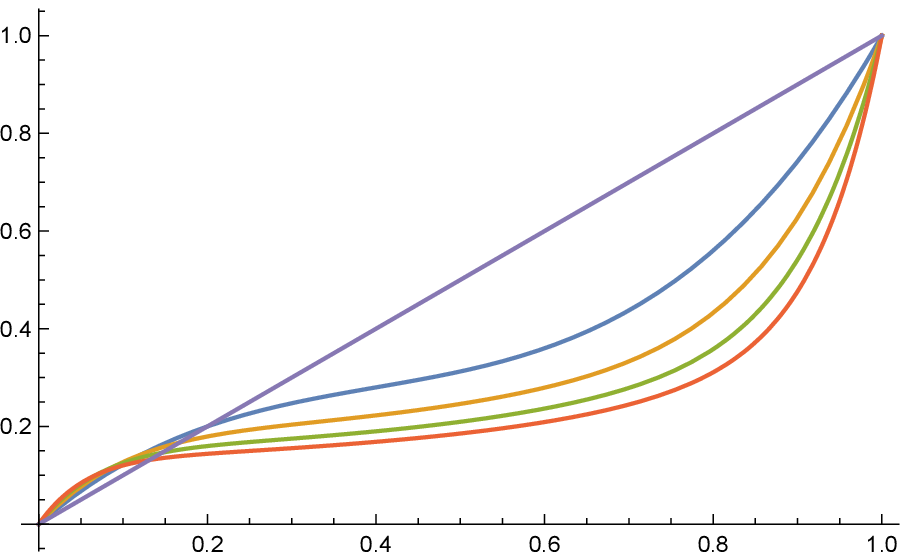}\includegraphics[width=7cm]{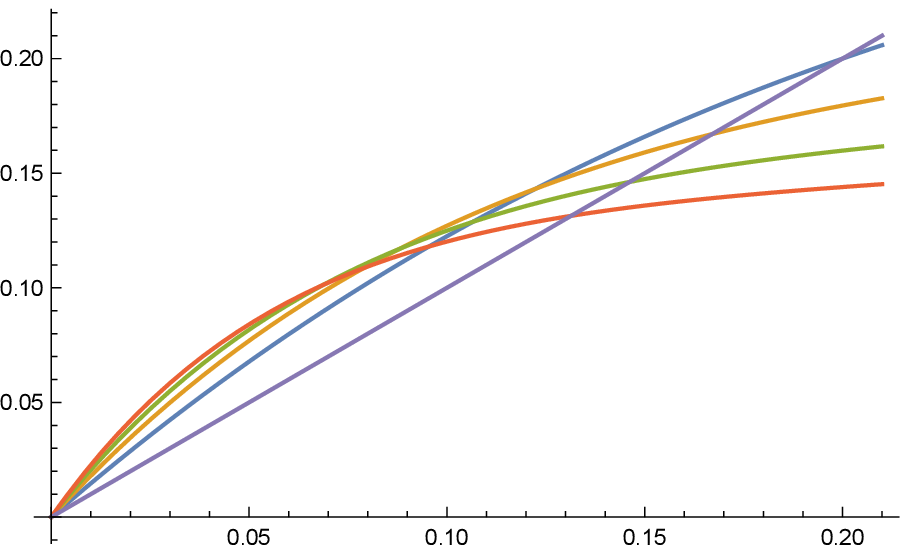}
\end{center}
\caption{The graphs of $f_n$ for $n=3$ (blue), $n=5$ (orange), $n=7$ (green) and $n=9$ (red),
on~$[0,1]$ and $[0,\nicefrac{1}{5}]$}
\end{figure}

\begin{Pro}
For all odd $n\ge 2$, $(0,1)$ is the basin of attraction  of $\alpha_n$.
\end{Pro}

\begin{proof}
The unicity of the fixed point $\alpha_n$ and formulas \eqref{obvious}-\eqref{obvious3} imply 
\begin{equation}
\label{eq.f-odd}
f_n(x)> x, \forall x\in (0,\alpha_n) \ \ \mbox{ and } \ \  f_n(x)< x, \forall x\in (\alpha_n, 1).
\end{equation}
Now we define the recursive sequence $(u_m)$ by $u_0=x_0$ and $u_{m+1}=f_n(u_m)$ for $m\ge0$. Since
\[\forall t\in[0,1]\,,\ \  f_n^\prime(t)=\frac{n}{\pi}\int_0^\pi(1-\cos x)(t(1-\cos x)+\cos x)^{n-1}\mathrm{d}x\ >\ 0,\]
the function $f_n$ is strictly increasing on $[0,1]$ and a simple reasoning shows that if $x_0\in(0,\alpha_n)$, $(u_m)$ is strictly increasing and bounded above by $\alpha_{n}$, and if $x_0\in(\alpha_{n},1)$, $(u_m)$ is strictly decreasing and bounded below by $\alpha_n$. As a consequence, for all $x_0\in(0,1)$:
\[\lim_{m\rightarrow\infty} u_m=\alpha_n. \]
\end{proof}

\begin{Rem}
\begin{enumerate}
\item Note that the fixed point $\alpha_{n}$ of $f_{n}$ is a linearly attracting, i.e. $f'_{n}(\alpha_{n})\in(-1,1)$.
Indeed, it holds obviously $f'_{n}(\alpha_{n})\in(0,1]$ by the preceding proof.
Moreover, the equality $f'_{n}(\alpha_{n})=1$ would imply that $\alpha_{n}$ is an inflection point of $f_{n}$ and then that $f^{(2)}_{n}(\alpha_{n})=0$, which would lead to $f_{n}(t)>t$ on $(\alpha_{n},1]$ since  $f^{(2)}_{n}$ is strictly increasing on $[0,1]$,
a contradiction.
\item The reasoning here can be applied for a GW with a reproduction law whose support is a subset of $2\mathbb N+1$. Indeed, the studied function $f=\sum_{n\geq 1}q_{2n+1}f_{2n+1}$
is then  strictly increasing on~$[0,1]$ and admits a unique fixed point on this interval. 
\end{enumerate}
\end{Rem}

\subsubsection{Proof of Theorem~\ref{only} and Corollary \ref{Corcool}}

\noindent 
{\bf The case of a {\boldmath$n$}-ary tree when {\boldmath$n\geq 3$} is odd}\medskip

Note first that in this case, the statement of Theorem~\ref{only} is obvious when 
 $\bp_{1}=\bp_{2}=\frac12$, since~$0$ is a fixed point  of $f_{n}$
and thus $(0,\frac12,\frac12,0_{k-2})$ is a fixed point of $H$.

It thus remains to prove Theorem~\ref{only} in this case when 
$\bp_{1}=\bp_{2}\in(0,\frac12)$. To this end,
let us fix $\bp \in \mathscr Q_{k}$ (with $k\geq 2$)
such that $\bp_{1}=\bp_{2}\in(0,\frac12)$ and $\bp_{2}>\bp_{3}\geq\dots\geq\bp_{k}\geq0$, 
and let us assume for a moment that 
there exists $\beta'>0$ such that $\bp_{0}(m)\geq \beta'$ for every $m>0$.
It then holds $0<\beta\leq \bp_{1}(m)\leq \frac{1-\beta'}2<\frac12$ for every $m>0$.
Thus, with the same arguments as those used in the proof of Proposition~\ref{nonunif_CV},
but working now with the compact set $[\beta,\frac{1-\beta'}2]\subset(0,\frac12)$ instead of $[\beta,\frac12]$,
one shows that 
$$\bp(m)\longrightarrow \left(\alpha_{n},\frac{1-\alpha_{n}}2, \frac{1-\alpha_{n}}2 , 0_{k-2}\right)\ \  \text{when $m\to \infty$}.$$

To conclude, let us then prove that
when $\bp_{1}=\bp_{2}\in(0,\frac12)$, there exists $\beta'>0$ such that $\bp_{0}(m)\geq \beta'$ for every $m>0$.

First, let us observe from the spreading rules that if $\bp_{0}(\ell)>0$ for some~$\ell\in\mathbb N$, then
$\bp_{0}(m)>0$ for every $m\geq\ell$. In particular, $\bp_{0}(m)>0$ for every $m\in \mathbb N$
when $\bp_{0}>0$ and, when $\bp_{0}=0$, then $k\geq 3$ and $\bp_{3}>0$, which implies $\bp_{0}(1)>0$ (also from the spreading rules, since $n\geq3$ is odd).
Consequently: $\bp_{0}(m)>0$ for every $m>0$. 

Moreover, note from the spreading rules that for every $m>0$,
\begin{align}
\label{eq.comp-fn}
\bp_{0}(m+1)=H_0(\bp_{0}(m),\dots,\bp_{k}(m))&\geq 
\sum_{k,0\le 2k\le n}\binom{n}{2k}\binom{2k}{k}\bp_{1}^{2k}(m)\bp_{0}^{n-2k}(m)\,.
\end{align}
Using now Remark~\ref{re.expo0} and $1-\bp_{0}(m)\geq 2\beta$, note also that
there exist $C=2\beta D>0$ and $a\in(0,1)$ such that for every $m>0$,
$\sum_{\ell=3}^{k}\bp_{\ell}(m)\leq Ca^{m}$ and thus
$$2\bp_{1}(m)=1-\bp_{0}(m)-\sum_{\ell=3}^{k}\bp_{\ell}(m)\geq
 \big(1-\bp_{0}(m)\big)\big(1-\frac{1}{2\beta}\sum_{\ell=3}^{k}\bp_{\ell}(m)\big)\geq
 \big(1-\bp_{0}(m)\big)\big(1 - Da^{m} \big)\,.$$
Take $m_{0}\in\mathbb N^{*}$ and $b\in(a,1)$ such that $Da^{m}\leq b^{m}$ for every $m\geq m_{0}$.
It then follows  from \eqref{eq.f-odd} and \eqref{eq.comp-fn} that:
\begin{align*}
\forall m\geq m_{0}\,,\ \ \ \bp_{0}(m+1)&\ \geq\  (1-b^{m})^{n}\sum_{k,0\le 2k\le n}\binom{n}{2k}\binom{2k}{k}\left(\frac{1-\bp_{0}(m)}2\right)^{2k}\bp_{0}^{n-2k}(m)\\
&\ =\  (1-b^{m})^{n}f_{n}(\bp_{0}(m))\ \geq\ (1-b^{m})^{n}\min\{\bp_{0}(m),\alpha_{n}\}.
\end{align*}
Reasoning by induction thus leads to:
\begin{align*}
\forall m\geq m_{0}\,,\ \ \ \bp_{0}(m)\ \geq\ \min\{\bp_{0}(m_{0}),\alpha_{n}\}\prod_{\ell=m_{0}}^{m-1}(1-b^{\ell})^{n}\ \geq \ \min\{\bp_{0}(m_{0}),\alpha_{n}\}\,A^{n}\,,
\end{align*}
where 
$A := \prod_{\ell=m_{0}}^{+\infty}(1-b^{\ell})$ is positive since 
the convergence of the Neumann series $\sum_{\ell\geq m_{0}}b^{\ell}$
implies the one of $\sum_{\ell\geq m_{0}}\ln(1-b^{\ell})$  to the real negative number $B=\ln(A)$.
It follows that for every $m>0$, 
$$\bp_{0}(m)\ \geq\  \beta':=\min\{A^{n}\alpha_{n},A^{n}\bp_{0}(m_{0}),\bp_{0}(m_{0}-1),\dots,\bp_{0}(1) \}\ >\ 0\,,$$ which concludes the proof of Theorem~\ref{only} 
in the case of a $n$-ary tree when $n\geq 3$ is odd.\\

\noindent 
{\bf The general case of a GW tree supported in {\boldmath$2\mathbb N+1$}}\medskip

We now look at the function $f=\sum_{n\geq 1}q_{2n+1}f_{2n+1}$
and at the corresponding function~$H$.
As above, the statement of Corollary~\ref{Corcool}  is obvious when 
 $\bp_{1}=\bp_{2}=\frac12$, since~$0$ is a fixed point  of $f$
and thus $(0,\frac12,\frac12,0_{k-2})$ is a fixed point of $H$.

It thus just remains to prove it when
$\bp_{1}=\bp_{2}\in(0,\frac12)$, so we fix $\bp \in \mathscr Q_{k}$ (with $k\geq 2$)
such that $\bp_{1}=\bp_{2}\in(0,\frac12)$ and $\bp_{2}>\bp_{3}\geq\dots\geq\bp_{k}\geq0$.
Reasoning as we did above with a $n$-ary tree when~$n\geq 3$ is odd,
it is sufficient to show  
 that 
there exists $\beta'>0$ such that $\bp_{0}(m)\geq \beta'$ for every~$m>0$.

To this end, note first  that the relation $\sum\limits_{n\geq 1}q_{2n+1}f'_{2n+1}(0)>1$ (see \eqref{obvious2})
implies
the existence 
  of $n^{*}\in\mathbb N^{*}$
such that $\sum\limits_{1\leq n\leq n^{*}}q_{2n+1}f'_{2n+1}(0)>1$.
The function $\tilde f:=\sum\limits_{1\leq n\leq n^{*}}q_{2n+1}f_{2n+1}$
hence satisfies $\tilde f(0)=0$, $\tilde f'(0)>1$, and $\tilde f\leq f$ on $[0,1]$, implying
$\tilde f(1)\leq f(1)=1 $.
It thus admits at least one fixed point in $(0,1)$
and we define $\alpha^{*}$ as the smallest one. It follows that
$\tilde f(x)>x$ on $(0,\alpha^{*})$ and, since $\tilde f$ is increasing on $[0,1]$,
the function $\tilde f$ satisfies $\tilde f(x)\geq \min\{x,\alpha^{*}\}$ for every $x\in[0,1]$.

We can then conclude by following the same lines as above for a $n$-ary tree:
again, the spreading rules imply that $\bp_{0}(m)>0$ for every $m>0$
and that
\begin{align*}
\bp_{0}(m+1)=H_0(\bp_{0}(m),\dots,\bp_{k}(m))&\geq 
\sum_{n\geq 1}q_{2n+1}\sum_{k,0\le 2k\le 2n+1}\binom{2n+1}{2k}\binom{2k}{k}\bp_{1}^{2k}(m)\bp_{0}^{2n+1-2k}(m)\\
&\geq 
\sum_{ n=1}^{n^{*}}q_{2n+1}\sum_{k,0\le 2k\le 2n+1}\binom{2n+1}{2k}\binom{2k}{k}\bp_{1}^{2k}(m)\bp_{0}^{2n+1-2k}(m)
\,.
\end{align*}
Reasoning as in the lines following \eqref{eq.comp-fn}
then implies the existence of  $m_{0}\in\mathbb N^{*}$
and of
 $b\in(0,1)$ such that, for every $m\geq m_{0}$,
\begin{align*}
\bp_{0}(m+1)&\ \geq\  
\sum_{ n=1}^{n^{*}}q_{2n+1}(1-b^{m})^{2n+1}\sum_{k,0\le 2k\le 2n+1}\binom{2n+1}{2k}\binom{2k}{k}\left(\frac{1-\bp_{0}(m)}2\right)^{2k}\bp_{0}^{2n+1-2k}(m)\\
&\ \geq\  (1-b^{m})^{2n^{*}+1}\tilde f(\bp_{0}(m))\ \geq\ (1-b^{m})^{2n^{*}+1}\min\{\bp_{0}(m),\alpha^{*}\}
\end{align*}
and then
\begin{align*}
\bp_{0}(m)\ \geq\ \min\{\bp_{0}(m_{0}),\alpha^{*}\}\prod_{\ell=m_{0}}^{m-1}(1-b^{\ell})^{2n^{*}+1}\ \geq \ \min\{\bp_{0}(m_{0}),\alpha^{*}\}\,\prod_{\ell=m_{0}}^{+\infty}(1-b^{\ell})^{2n^{*}+1}\ >\ 0\,.
\end{align*}
This implies the existence of $\beta'>0$ such that 
$\bp_{0}(m)\geq \beta'$ for every $m>0$ and then concludes the proof of Corollary \ref{Corcool}.


\subsection{The even case}

\subsubsection{The fixed point  {\boldmath$\alpha_{n}$} is  linearly attracting} 
\label{sub.even-attract}

\begin{figure}[h!]
\begin{center}
\includegraphics[width=7cm]{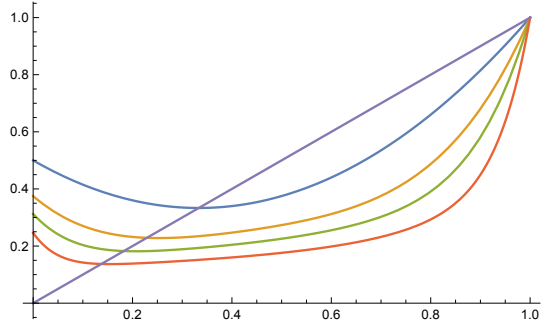}
\end{center}
\caption{The graphs of $f_n$ for $n=2$ (blue), $n=4$ (orange), $n=6$ (green) and $n=10$ (red)}
\end{figure}

\begin{Pro}\label{evenattractive}
For all even $n\ge 2$, $\alpha_n$ is linearly attracting. 
\end{Pro}

The only difficulty to obtain this statement is to prove that $f^\prime_n(\alpha_n)>-1$. Indeed, since $f_n(0)>0$ and $f_n'(1)>1=f_n(1)$, the unicity of $\alpha_{n }$ 
leads to
\[f_n(x)> x, \forall x\in (0,\alpha_n)\ \  \mbox{ and }  \ \ f_n(x)< x, \forall x\in (\alpha_n, 1)\,,\]
and hence $f'(\alpha_n)\le 1$.
Since moreover $f$ is strictly convex on $[0,1]$, we have $f'(\alpha_n)< 1$ since the equality 
$f'(\alpha_n)=1$ would imply $f_n(x)> x$ on $[0,1]$, a contradiction.


A direct proof of Proposition \ref{evenattractive} using integral estimates
relying on the relation \eqref{Luc} is proposed in the following subsection. Nevertheless, we would like to point out that we have come up with a totally independent proof using Budan's theorem:

\begin{The}{\bf(of Budan-Fourier)}\cite{Budan}\label{budan}
Let $P(x)=0$ be a polynomial equation with real coefficients  of degree $n$ and let $a<b$ be any two real numbers. Then,
there exists $k\in\N$ such that
 the number of roots (counted with multiplicity) of this equation in the interval $(a,b]$ is equal to $$V_a(P)-V_b(P)-2k\,,$$ where, 
for $c\in\R$, $V_c(P)$ is the number of sign variations  in the sequence $P(c),P^\prime(c),\dots,P^{(n)}(c)$.
\end{The}



This proof has its own interest since it can be applied to prove the attractivity of a fixed point in the more general setting of GW, see Remark \ref{re.GW-n}. It relies on the



\begin{Lem}
\label{le.Budan}
For all even $n\ge 2$, the function
\[\gamma:t\in[0,1]\longmapsto tf_n(t)=\sum_{k,0\le 2k\le n}\binom{n}{2k}\binom{2k}{k}\left(\frac{1-t}{2}\right)^{2k}t^{n+1-2k}\in\mathbb R,\]
is strictly increasing on $(0,\nicefrac{1}{2})$.
\end{Lem}

\begin{proof}
Writing:
\begin{align*}
\left(\left({1-t}\right)^{2k}t^{n+1-2k}\right)^\prime
&=t^n\left(\frac{1-t}{t}\right)^{2k-1}\left(\frac{1}{t}(n+1-2k)-(n+1)\right),
\end{align*}
the inequality $\gamma^\prime(t)>0$ for $t\in(0,\nicefrac12)$ is equivalent to: 
\[\forall\,t\,\in\,\left(0,\nicefrac12\right)\,,\ \ \ \sum_{k,0\le 2k\le n}2^{-2k}\binom{n}{2k}\binom{2k}{k}\left(\frac{1-t}{t}\right)^{2k-1}\left(\frac{1}{t}(n+1-2k)-(n+1)\right)>0.\]
Using the substitution $s=\frac{1-t}{t}\Leftrightarrow t=\frac{1}{1+s}$, it is equivalent to prove on $(1,+\infty)$:
\begin{align*}
g(s):&=\sum_{k,0\le 2k\le n}2^{-2k}\binom{n}{2k}\binom{2k}{k}s^{2k-1}\left((1+s)(n+1-2k)-(n+1)\right)\\
&=\sum_{k,0\le 2k\le n}2^{-2k}\binom{n}{2k}\binom{2k}{k}s^{2k-1}\left(s(n+1-2k)-2k\right)>0.
\end{align*}

In order to use Theorem \ref{budan}, we need the $\ell$-th derivatives of the function $g$ for $0\le \ell\le n$ and their values at the limits of the interval $(1,+\infty)$. 
As:
\begin{align*}
\frac{\mathrm{d}^\ell}{\mathrm{d}s^\ell}(s^{2k})&=\frac{(2k)!}{(2k-\ell)!}s^{2k-\ell},\\
\frac{\mathrm{d}^\ell}{\mathrm{d}s^\ell}(2ks^{2k-1})&=\frac{(2k)!}{(2k-\ell)!}(2k-\ell)s^{2k-\ell-1},\\
\binom{n}{2k}&=\frac{n!}{(n-\ell)!}\binom{n-\ell}{2k-\ell}\frac{(2k-\ell)!}{(2k)!},
\end{align*}
we obtain:
\begin{align}\label{lieme}
g^{(\ell)}(s)
=\frac{n!}{(n-\ell)!}\sum_{k,\ell\le 2k\le n}2^{-2k}\binom{2k}{k}\binom{n-\ell}{2k-\ell}s^{2k-\ell-1}\left(s(n+1-2k)-(2k-\ell)\right).
\end{align}
Since $n$ is even and $2k\leq  n$ for every $k$  considered in the sum in \eqref{lieme},
\[g^{(\ell)}(s)\underset{s\rightarrow\infty}{\sim}s^{n-\ell}2^{-n}\binom{n}{\frac n2}\frac{n!}{(n-\ell)!}>0, \forall \ell \in\llbracket 0,n\rrbracket .\]
Moreover, according again to \eqref{lieme}, for every $\ell\in\llbracket 0,n\rrbracket$:
\[g^{(\ell)}(1)=\frac{n!}{(n-\ell)!}\sum_{k,\ell\le 2k\le n}2^{-2k}\binom{2k}{k}\binom{n-\ell}{2k-\ell}\left(n+1+\ell-4k\right). \]
The sign of
$g^{(\ell)}(1)$, and thus of
\[\sum_{k,\ell\le 2k\le n}2^{-2k}\binom{2k}{k}\binom{n-\ell}{2k-\ell}\left(n+1+\ell-4k\right):=\sum_{k,\ell\le 2k\le n}\mu_k\alpha_k\,,\]
where $\alpha_k:=n+1+\ell-4k$, 
is difficult to obtain directly. As $(\alpha_k)_{k\geq0}$ is a decreasing sequence, the main idea is to use Lemma \ref{wiard1} to bound below this sum with a quantity that we are able to compute. 
Defining
\[\nu_k:=\binom{n-\ell}{2k-\ell},\ \  \text{it holds}\ \  \frac{\mu_k}{\nu_k}={2}^{-2k}\binom{2k}{k}=\xi_{2k}\,, \]
and, as $\left(\xi_{2k}\right)_{k\ge0}$ is a decreasing sequence (see Lemma \ref{XI}),  
we can apply Lemma \ref{wiard1} which gives:
\begin{align*}
\frac{\sum_{k,\ell\le 2k\le n}\mu_k\alpha_k}{\sum_{k,\ell\le 2k\le n}\mu_k}
\ge \frac{\sum_{k,\ell\le 2k\le n}\nu_k\alpha_k}{\sum_{k,\ell\le 2k\le n}\nu_k}=\frac{\sum_{k,\ell\le 2k\le n}\binom{n-\ell}{2k-\ell}(n+1+\ell-4k)}{\sum_{k,\ell\le 2k\le n}\binom{n-\ell}{2k-\ell}}.
\end{align*}
Moreover, 
according to  \eqref{combineasyimpli} and to \eqref{eq.binom-easy}:
\begin{align*}
\sum_{k,\ell\le 2k\le n}\nu_k\alpha_k&=(n-\ell+1)\sum_{k,\ell\le 2k\le n}\binom{n-\ell}{2k-\ell}-2\sum_{k,\ell\le 2k\le n}\binom{n-\ell}{2k-\ell}(2k-\ell)\\
&=
\begin{cases}(n-\ell+1)2^{n-\ell-1}-2(n-\ell)2^{n-\ell-2}\ &\text{if}\ \ell\in\{0,\dots,n-2\},\\
(n-\ell+1)2^{n-\ell-1}-2\ &\text{if}\ \ell=n-1,\\
1\ &\text{if}\ \ell=n,
\end{cases}\\
&=
\begin{cases}2^{n-\ell-1}\ &\text{if}\ \ell\in\{0,\dots,n-2\},\\
0\ &\text{if}\ \ell=n-1,\\
1\ &\text{if}\ \ell=n.
\end{cases}
\end{align*}
It follows that $g^{(\ell)}(1)\ge0$ for every $\ell\in\llbracket 0,n\rrbracket$,
so
we have proved that $g^{(\ell)}(1)$ and $g^{(\ell)}(+\infty)$ have always the same sign. According to Theorem~\ref{budan}, the number of roots of $g$ in $(1,+\infty)$ is thus zero and hence $g>0$ on $(1,+\infty)$ since $\lim\limits_{t\to+\infty} g(t)=+\infty$.
\end{proof}

\

\noindent
{\bf Proof of Proposition \ref{evenattractive}}\\
Using $f_{n}(\alpha_{n})<1$ and Lemma~\ref{le.Budan}, the proof of  Proposition \ref{evenattractive} is straightforward:
\[\forall\,t\,\in\,(0,\nicefrac{1}{2})\,,\ \ \ \gamma^\prime(t)=f_n(t)+tf_n^\prime(t)>0 \Leftrightarrow f_n^\prime(t)> - \frac{f_n(t)}{t},\]
and taking $t=\alpha_{n}$ leads to $f_n^\prime(\alpha_n)>-1$.

\begin{Rem}
\label{re.GW-n}

The statement of Lemma~\ref{le.Budan} remains actually true when $n\geq2$ is odd.
It follows that
in the GW case, the function 
\[t\mapsto tf(t)=\sum_{n\geq 2}q_n t f_n(t)\]
is strictly increasing on $(0,\nicefrac{1}{2})$. 
In particular, we have $f^\prime(\alpha)>-1$ for every fixed point $\alpha\in(0,\frac12)$.\\ 
Recall moreover (Remark~\ref{re.GW1}) that
$f$ admits at least one fixed point in $(0,\frac12)$ and 
that, either $f(0)=0$ and $f'(0)>1$, or $f(0)>0$. 
Hence, denoting by $\alpha$ the smallest fixed point in $(0,\frac12)$, 
we have necessarily $-1<f'(\alpha)\le 1$, which almost implies
the linear attractivity of $\alpha$. \\
Furthermore,  when the support of $q$ is included in $2\N$, the convexity of $f$ implies the attractivity of its unique fixed point.

\end{Rem}

\subsubsection{Basin of attraction of the fixed point \boldmath$\alpha_{n}$}

The attractivity of $\alpha_n$ is not enough to obtain the even case in Theorem~\ref{only}. 
In order to apply Proposition~\ref{nonunif_CV},
we have to prove that  the basin of attraction of $\alpha_n$ is $[0,1)$.

The proof is carried out in two steps. First, we prove
the existence of $n_{0}\in\mathbb N$ such
 that $f_{n}'(\alpha_n)>0$ for every even $n> n_{0}$,  which implies that the  basin of attraction of $\alpha_n$ is $[0,1)$
 when $n > n_{0}$. Secondly, we prove numerically that
the basin of attraction of $\alpha_n$
is also $[0,1)$ for every even $2\leq n \leq n_{0}$. Moreover, we  estimate the constants appearing in the computations with precision in order to minimize $n_{0}$ and then the number of values of $n$ for which we have to check the result numerically.

\begin{Lem}\label{DLP} 
For all  even $n\in \N^*$ and all $t\in[0,1]$,
\begin{equation}
f_n(t)\ge \frac{1}{\sqrt{2\pi (n+1)}}\,.
\end{equation}
\end{Lem}
\begin{proof}
Let us simply note that for every even $n\in \N^*$ and every $t\in[0,1]$,
\begin{eqnarray*}
f_n(t)
&=&\frac1{\pi} \int_0^{ \frac{\pi}{2}}\left((1-t)\cos{x}+t\right)^n{\rm d}x+\frac1{\pi} \int_{\frac{\pi}{2}}^{\pi}\left((1-t)\cos{x}+t\right)^n{\rm d}x\\
&\ge&\frac1{\pi} \int_0^{\frac{\pi}{2}}\cos^n{x}{\rm d}x\ge \frac{1}{\sqrt{2\pi (n+1)}},
\end{eqnarray*}
where the last inequality follows from standard estimates on Wallis integrals, see \eqref{easystuff}.
\end{proof}

\begin{Rem}
Lemma \ref{DLP} implies in particular $\alpha_n\ge  \frac{1}{\sqrt{2\pi (n+1)}}$ for every even $n\in\N^*$. 
\end{Rem}
Let us now recall  two classical results, obtained by integration by parts,
which will be useful in the sequel:
if $X$  follows the standard normal distribution, then
\begin{align}
\forall x>0\,,\ \ \p(X\ge x)&=\frac{1}{\sqrt{2\pi}}\int_x^\infty e^{-\frac{u^2}{2}}\mathrm{d}u\le \frac{e^{-\frac{x^2}{2}}}{x\sqrt{2\pi}},\label{foncrep} \\
\text{and}\qquad \forall n\in\mathbb N\,,\ \  \E\left[X^{2n}\right]&=\frac{(2n)!}{2^n n!}.\label{moment}
\end{align} 
On one hand, using \eqref{foncrep}:
\begin{align*}
\frac1{\pi} \int_0^{ 1}e^{-n\frac{x^2}{2}}{\rm d}x&=\sqrt{\frac{2}{n\pi}}\left(\frac{1}{\sqrt{2\pi}}\int_0^{\sqrt n}e^{-\frac{u^2}{2}}\mathrm{d}u\right)=\sqrt{\frac{2}{n\pi}}\left(\frac{1}{2}-\p(X\ge \sqrt{n})\right)\\
&\ge \sqrt{\frac{2}{n\pi}}\left(\frac{1}{2}-\frac{e^{-\frac{n}{2}}}{\sqrt{2\pi n}}\right)=\frac{1}{\sqrt{2\pi n}}-\frac{e^{-\frac{n}{2}}}{\pi n}.
\end{align*}

\begin{Lem}\label{DPA}
There exists $n_0\in\mathbb N$ such that for all even $n>n_0$, we have $f_n'\left( \frac{1}{\sqrt{2\pi (n+1)}}\right)>0$.
\end{Lem}

\begin{proof}
Recall that:
\[f_n'(x)= \frac{n}{\pi} \int_0^{\pi}(1-\cos{t})\left((1-x)\cos{t}+x\right)^{n-1}{\rm d}t.\]
Let $\zeta_n= \frac{1}{\sqrt{2 \pi (n+1)}}$. 
We have
\begin{eqnarray*}
\frac{\pi}{n}f_n'(\zeta_n)=& \displaystyle\int_0^{1}(1-\cos{t})\left((1-\zeta_n)\cos{t}+\zeta_n\right)^{n-1}{\rm d}t&+\int_1^{\pi}(1-\cos{t})\left((1-\zeta_n)\cos{t}+\zeta_n\right)^{n-1}{\rm d}t\\
=&A&+\hspace{3cm}B.
\end{eqnarray*}
Using successively that  
$\cos u\le 1-\nicefrac{u^2}{2}+\nicefrac{u^4}{4!}$,  $\cos u\ge 1-\nicefrac{u^2}{2}$, and
$\ln (1-\nicefrac{u^2}{2})+\nicefrac{u^2}{2}\ge- \nicefrac{u^4}{4}$
on $[0,1]$, and $e^u-1\ge u$ on $\mathbb R$,
we obtain
\begin{eqnarray*}
A&=&\int_0^{1}(1-\cos{t})\left((1-\zeta_n)\cos{t}+\zeta_n\right)^{n-1}{\rm d}t\ge \int_0^{1} \left(\frac{t^2}{2}-\frac{t^4}{4!}\right)e^{(n-1)\ln(\cos{t})}{\rm d}t\\
&\ge& \int_0^{1} \left(\frac{t^2}{2}-\frac{t^4}{4!}\right)e^{-(n-1)\frac{t^2}{2}}{\rm d}t-(n-1)\int_0^{ 1}  \left(\frac{t^2}{2}-\frac{t^4}{4!}\right)\frac{t^4}{4}  e^{-(n-1)\frac{t^2}{2}}  {\rm d}t\\
&\ge& \int_0^{1} \left(\frac{t^2}{2}-\frac{t^4}{4!}\right)e^{-(n-1)\frac{t^2}{2}}{\rm d}t-(n-1)\int_0^{ 1}  \frac{t^6}{8} e^{-(n-1)\frac{t^2}{2}}  {\rm d}t \\
&\ge& \int_0^{1} \frac{t^2}{2}e^{-(n-1)\frac{t^2}{2}}{\rm d}t- \int_0^{1} \frac{t^4}{4!}e^{-(n-1)\frac{t^2}{2}}{\rm d}t-(n-1)\int_0^{ 1}  \frac{t^6}{8} e^{-(n-1)\frac{t^2}{2}}  {\rm d}t=:D-(E+F).  
\end{eqnarray*}
It is an easy task to see that with the substitution $u=t\sqrt{n-1}$ and \eqref{moment}:
\begin{align}\label{EF}
\frac{E+F}{\sqrt{2\pi}}\le(n-1)^{-\nicefrac{5}{2}}\left(\frac{\E[X^4]}{2\times 4! }+\frac{\E[X^6]}{16}\right)=(n-1)^{-\nicefrac{5}{2}}.
\end{align}
In order to bound below $D$, we use integration by parts, \eqref{moment} and \eqref{foncrep}:
\begin{align}\label{A}
 \frac{D}{\sqrt{2\pi}}&= (n-1)^{-\nicefrac{3}{2}}\frac{1}{2\sqrt{2\pi}} \int_0^{\sqrt{n-1}}u^2e^{-\frac{u^2}{2}}{\rm d}u = (n-1)^{-\nicefrac{3}{2}}\left(\frac{\E[X^2]}{4}-\frac{1}{2\sqrt{2\pi}} \int_{\sqrt{n-1}}^\infty u^2e^{-\frac{u^2}{2}}{\rm d}u\right)\nonumber\\
 &=(n-1)^{-\nicefrac{3}{2}}\left(\frac{1}{4}-\frac{1}{2\sqrt{2\pi}} \left(\sqrt{n-1}e^{-\frac{n-1}{2}}+\int_{\sqrt{n-1}}^\infty e^{-\frac{u^2}{2}}{\rm d}u\right)\right)\nonumber\\
 &\ge (n-1)^{-\nicefrac{3}{2}}\left(\frac{1}{4}-\frac{1}{2\sqrt{2\pi}} \left(\sqrt{n-1}e^{-\frac{n-1}{2}}+ \frac{e^{-\frac{n-1}{2}}}{\sqrt{n-1}}\right)\right)\nonumber\\
 &= (n-1)^{-\nicefrac{3}{2}}\left(\frac{1}{4}-\frac{n}{2\sqrt{2\pi(n-1)}} e^{-\frac{n-1}{2}}\right).
 \end{align}

Let us now turn to $B$ and denote by $\mathscr M:=\left\{t\in[1,\pi], (1-\zeta_n)\cos{t}+\zeta_n<0\right\}$. Since $n-1$ is odd, we have:

\begin{align*}
B&= \int_1^{\pi}(1-\cos{t})\left((1-\zeta_n)\cos{t}+\zeta_n\right)^{n-1}{\rm d}t  \ge  \int_{\mathscr M}(1-\cos{t})\left((1-\zeta_n)\cos{t}+\zeta_n\right)^{n-1}{\rm d}t\\
&\ge \int_{\frac{\pi}{2}}^{\pi}(1-\cos{t}) (1-\zeta_n)^{n-1}\left(\cos{t}\right)^{n-1}{\rm d}t=- (1-\zeta_n)^{n-1}\left(W_n+W_{n-1}\right),
\end{align*}
 where $W_n$ is the Wallis integral. As $(1-\zeta_n)^{n-1}\le e^{-{\frac{n-1}{\sqrt{2\pi (n+1)}}} }$ and $(W_n)$ is a decreasing sequence such that for all $n\geq 1$, $W_{n}\le \sqrt{\frac{\pi}{2n}}$:
\begin{equation}\label{B}
B\ge -2W_{n-1}e^{-{\frac{n-1}{\sqrt{2\pi (n+1)}}} }\ge -\sqrt{\frac{{2\pi}}{n-1}}e^{-{\frac{n-1}{\sqrt{2\pi (n+1)}}} }.
\end{equation}

%
%

Combining \eqref{EF}, \eqref{A}, and \eqref{B} yields
\begin{align*}
\frac{\pi}{n}f_n'(\zeta_n)
&\ge \sqrt{2\pi} (n-1)^{-\nicefrac{3}{2}}\left(\frac{1}{4}-\frac{n}{2\sqrt{2\pi(n-1)}} e^{-\frac{n-1}{2}}-\frac{1}{(n-1)} -(n-1)e^{-{\frac{n-1}{\sqrt{2\pi (n+1)}}} }\right)\\
&=\sqrt{2\pi} (n-1)^{-\nicefrac{3}{2}}\left(\frac{1}{4}-w_n\right).
\end{align*}
The statement of Proposition~\ref{DPA} follows, since $\lim_{n\rightarrow\infty}w_n=0$. 
 \end{proof}
 
 \begin{figure}[h!]
\begin{center}
\includegraphics[width=7.5cm]{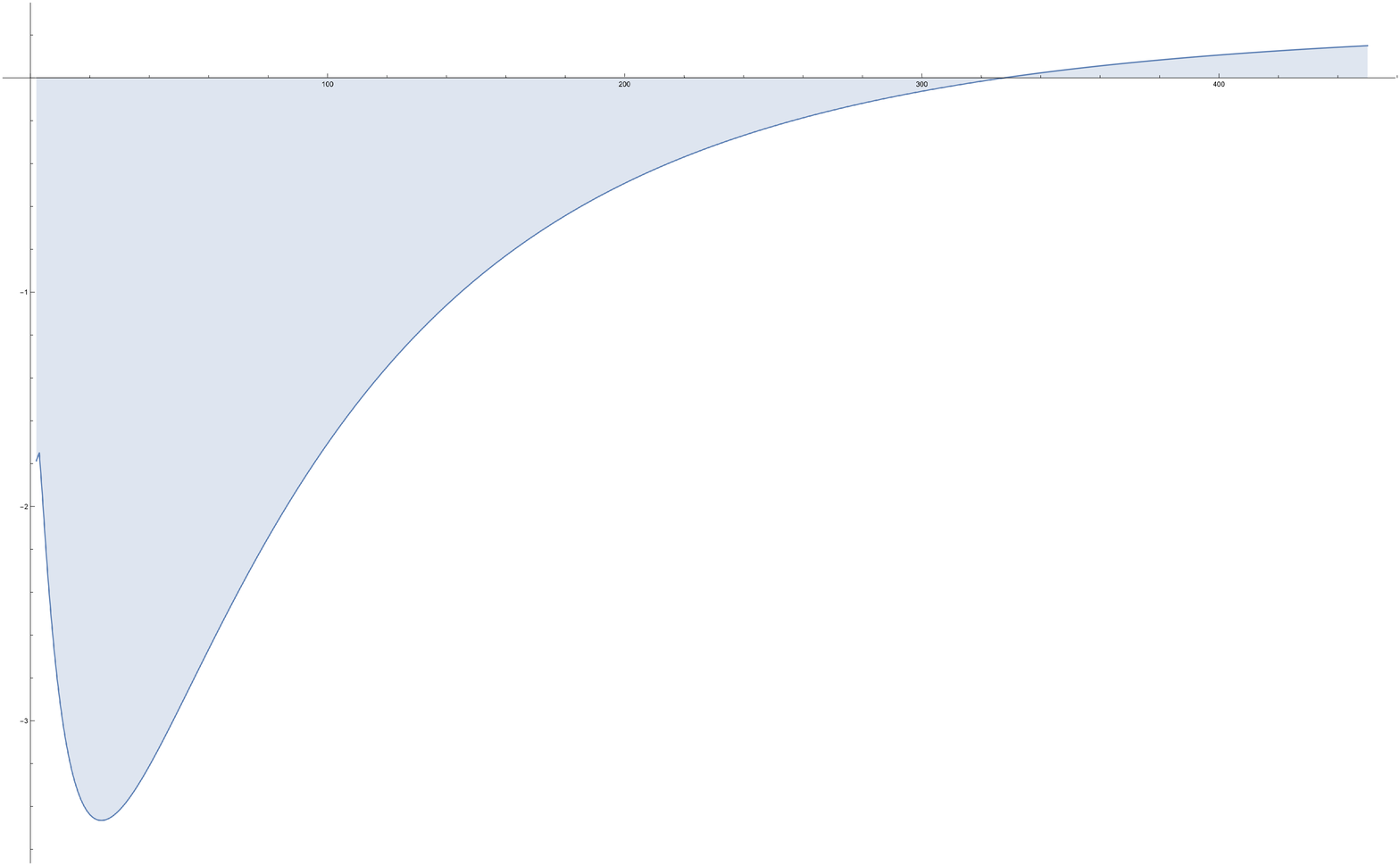}\includegraphics[width=7.5cm]{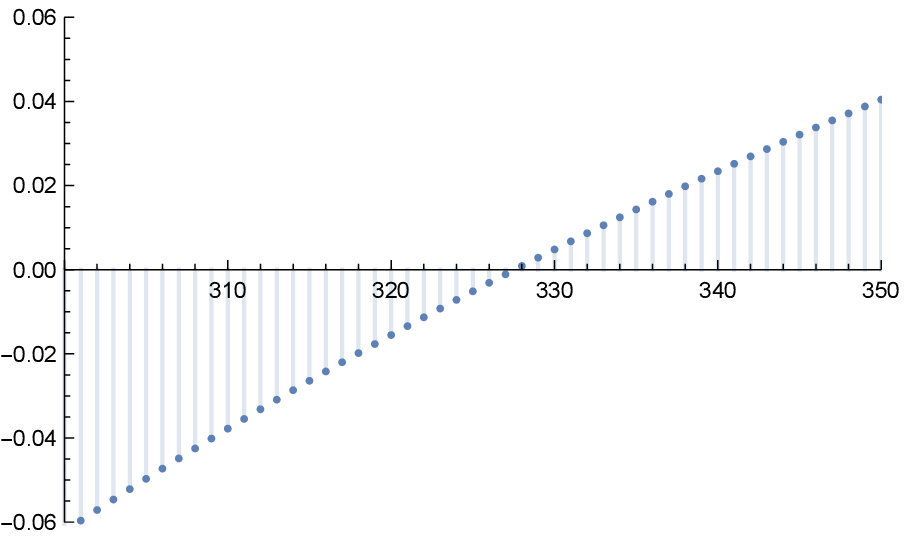}
\end{center}
\caption{$\frac{1}{4}-w_n$.}
\label{wn}
\end{figure}
 We introduce the following notation for every even $n\geq 2$:
 we denote by
\begin{equation}
\hat x_n:=\argmin_{x\in[0,1]}f_n(x)
\end{equation}
the global minimum of the function $f_n$ on $[0,1]$.
Note that the strict convexity of $f_n$ (as $n$ is even)
together with $f_{n}'(0)<0<f'_{n}(1)$ ensures that $\hat x_n$ is unique and belong to $(0,1)$, with
$f_{n}'(\hat x_n)=0$. 
 
\begin{Lem}\label{derivsup}
For every even  $n\in\mathbb N^{*}$ such that $f'_{n}(\alpha_{n})>0$,
the basin of attraction of $\alpha_n$ is $[0,1)$. 
This is in particular the case for every even $n$ sufficiently large. 
\end{Lem}
\begin{proof}
The second part of Lemma~\ref{derivsup}
is an immediate consequence  of its first part and of the convexity of
$f_n$ together with Lemmas \ref{DLP} and \ref{DPA}, which imply $f'_{n}(\alpha_{n})>0$
for every even $n$ large enough.

Let us then prove the first part of  this lemma.
Thanks to the convexity of $f_n$, we have $0<\hat x_n<\alpha_{n}$ and then $f_n(\hat x_n)> \hat x_n$.
We drop the subscript $_{n}$ in the rest of the proof below to lighten the notation
and we define the recursive sequence $(u_m)_{m\geq0}$ by $u_0=x_0\in[0,1)$ and $u_{m+1}=f_n(u_m)$ for $m\ge0$.\\
\underline{On  $[\alpha,1)$:} $f$ is increasing, $f([\alpha,1))=[\alpha,1)$, and $f(x)\le x$.\\ Thus, if $x_0\in[\alpha, 1)$, the sequence $(u_m)$ is decreasing and bounded below by $\alpha$,  implying that $(u_m)$ tends to $\alpha$, the fixed point of $f$ in~$[\alpha,1)$.\\
\underline{On $[f(\hat x),\alpha)\subset (\hat x,\alpha) $:}
 $f$ is increasing, $f([f(\hat x),\alpha))\subset [f(\hat x),\alpha)$, and $f(x)> x$.\\
Thus, if $x_0\in [f(\hat x),\alpha)$, the sequence $(u_m)$ is increasing and bounded above by $\alpha$,  
implying that $(u_m)$ tends to $\alpha$, the fixed point of $f$ in $ [f(\hat x),\alpha]$.\\
\underline{On  $[\hat x,f(\hat x))$:}  $f$ is increasing and $f([\hat x,f(\hat x))\subset [f(\hat x), 1)$.\\
We can thus conclude with the two previous cases when $x_{0}\in[\hat x,f(\hat x)) $.\\
\underline{On $[0,\hat x)$:} $f$ is decreasing and $f([0,\hat x))=(f(\hat x), f(0)]\subset [f(\hat x), 1)$.\\
We can thus again conclude with the two first cases when  $x_{0}\in[0,\hat x)$.
\end{proof}

\begin{Rem}\label{remark.1}
In the general GW setting, if the function $f=\sum_nq_nf_n$ is strictly convex, using previous arguments we get existence and unicity of the fixed point $\alpha$ in $[0,1)$. If in addition  $f^\prime(\alpha)\ge 0$,  the basin of attraction of $\alpha$ is $[0,1)$ with a similar reasoning as in the proof of Lemma \ref{derivsup}.
\end{Rem}
\subsubsection{Proof of Theorem~\ref{only} in the even case}

\noindent
{\bf Proof in the even case when  \boldmath$n>26$.}\medskip

Let us observe that the sequence $(w_n)$ converging to $0$ defined at  the end of the proof of Lemma~\ref{DPA}
is decreasing for $n$ large enough. More precisely, 

\begin{align*}
w_n&=\frac{\sqrt{n-1}}{2\sqrt{2\pi}} e^{-\frac{n-1}{2}}+\frac{1}{2\sqrt{2\pi}} e^{-\frac{n-1}{2}}+\frac{1}{(n-1)} +(n-1)e^{-{\frac{n-1}{\sqrt{2\pi (n+1)}}}}\\
&=w_n'+(n-1)e^{-{\frac{n-1}{\sqrt{2\pi (n+1)}}}}
\end{align*}
where $w_n'$ is decreasing for all $n>1$. \\
On the other hand, the derivative of the function $x\mapsto (x-1)e^{-{\frac{x-1}{\sqrt{2\pi (x+1)}}}}$ is decreasing for $x>\frac92\pi$ and negative at $x=26$ and hence $(w_n)$ is decreasing for $n\ge26$.

In particular, taking any  $n_{0}\geq 26$ such that  $\frac14-w_{n_{0}}>0$, 
one has $f_{n}'(\zeta_{n})>0$ and then $f_n'(\alpha_n)>0$ for every even $n\geq n_{0}$.
It is also easy to check numerically that $\frac14-w_{350}>0$ (see the right graph in Figure~\ref{wn}),
and thus  $f_n'(\alpha_n)>0$ for every even $n\ge 350$.

Moreover, computer assisted estimates show that  $f_n'(\alpha_n)>0$ for every even  $26<n<350$,
see Figure~\ref{fig-26} below.

Thus,  $f_n'(\alpha_n)>0$ for every even $n>26$ and  it follows from 
 Lemma~\ref{derivsup} that the basin of attraction of $\alpha_n$ is $[0,1)$
 for every even $n>26$. The statement of Theorem~\ref{only} in this case is then
 a consequence of Proposition~\ref{nonunif_CV}.
 
\begin{figure}[h!]
\begin{center}
\includegraphics[width=7cm]{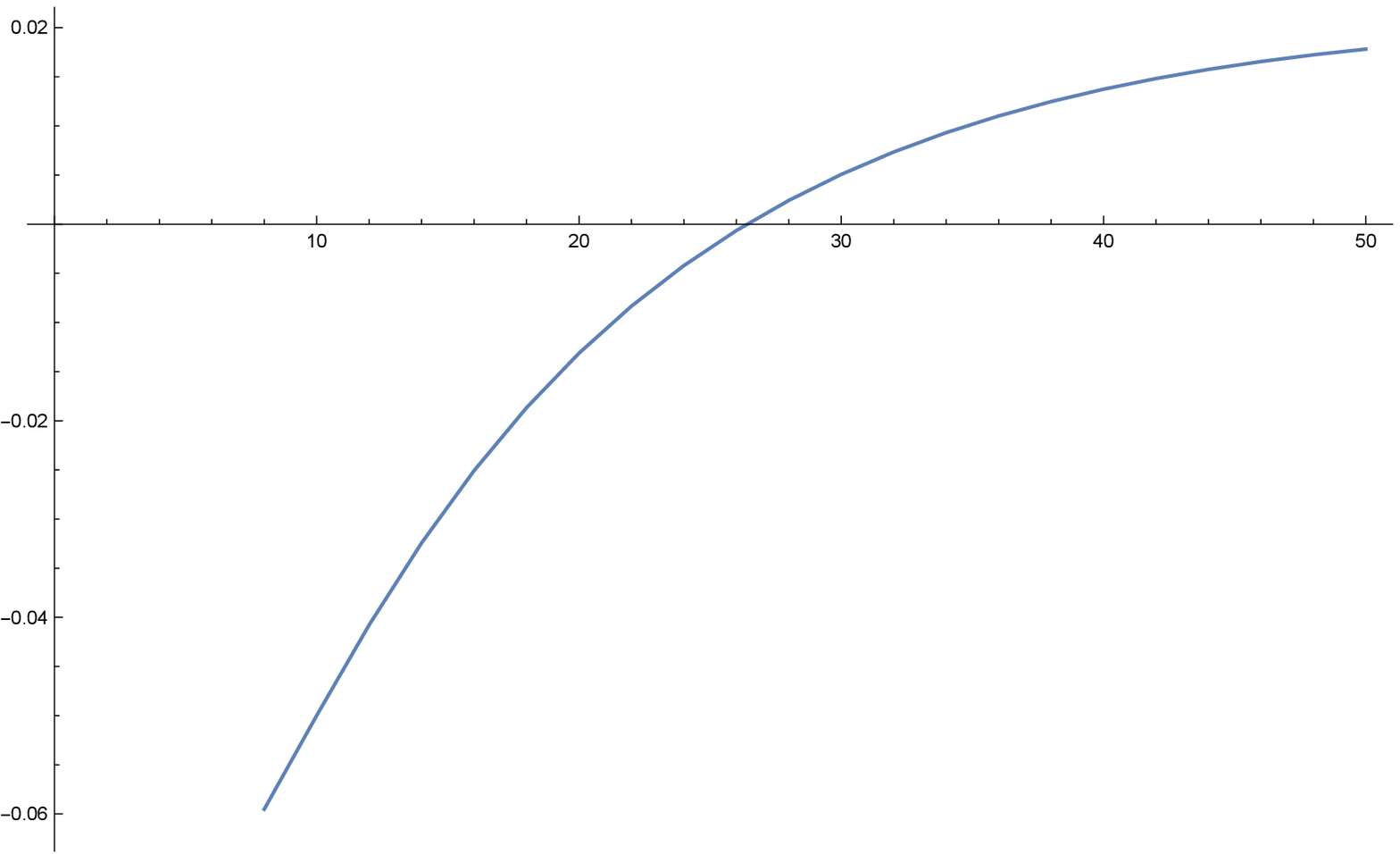}\includegraphics[width=7cm]{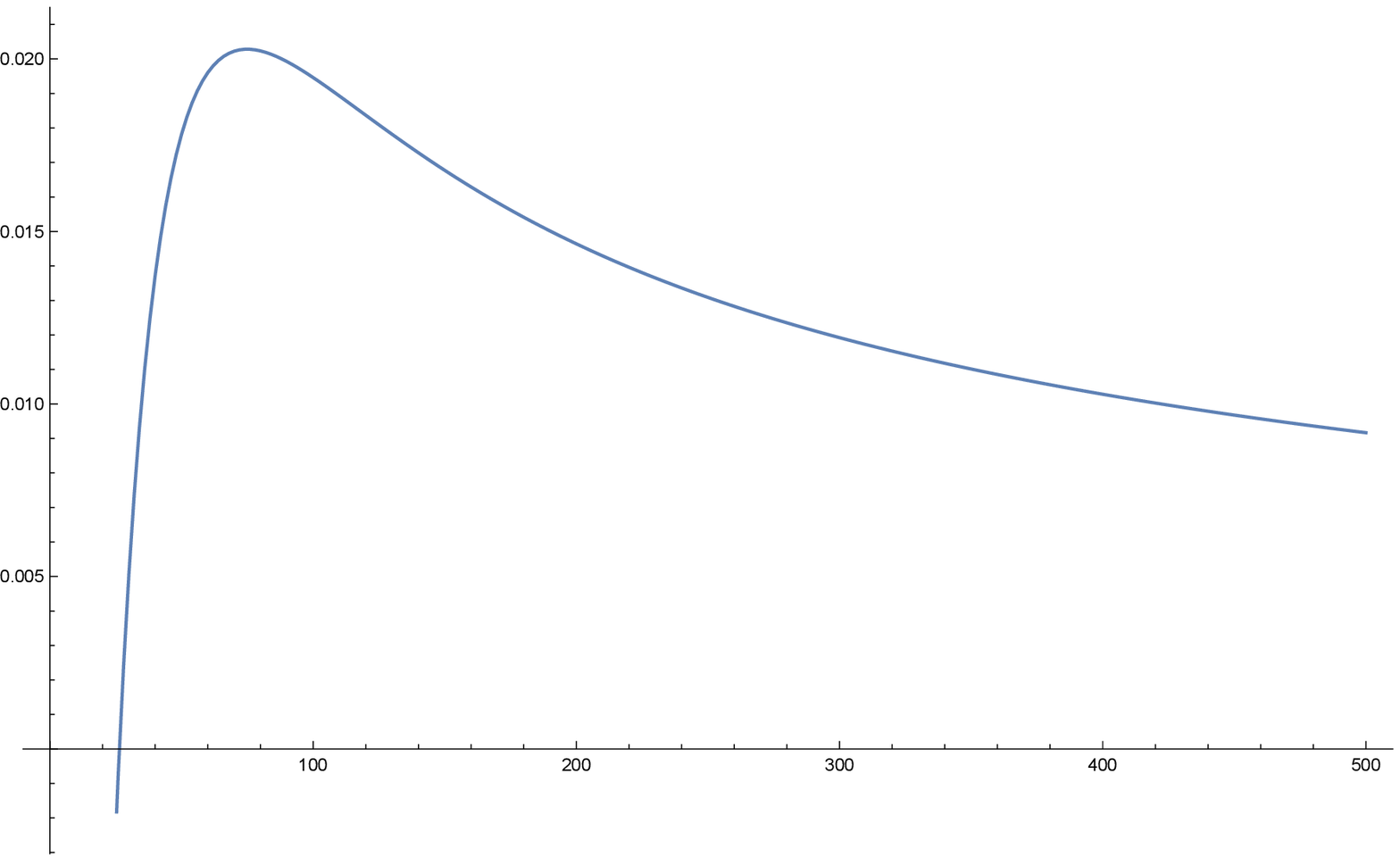}
\end{center}
\caption{$f^\prime_{2n}(\alpha_{2n})$ for $n$ in $[2,250]$}
\label{fig-26}
\end{figure}


\noindent
{\bf Proof of Theorem~\ref{only} in the even case when \boldmath$2\leq n\leq 26$.}\medskip

\begin{figure}[h!]
\begin{center}
\includegraphics[width=7cm]{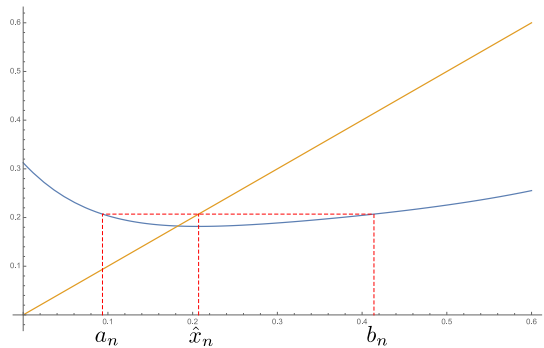}
\end{center}
\caption{The graph of $f_2$}
\label{f_2}
\end{figure}

    Let us consider the case  $n$ even and $f_{n}'(\alpha_n)<0$. The function $f_n$ being strictly convex  on $[0,1]$, the inverse image $f_n^{-1}(\hat{x}_n)$ of  its minimum $\hat{x}_n$ is composed by at most two elements, $a_n<b_n \in [0,1]$ (see Figure \ref{f_2}).

We have $\alpha_n<\hat{x}_n$ (since $f_{n}'(\alpha_n)<0$ and $f_{n}'(\hat{x}_n)=0$) and $f_{n}(\hat{x}_n)<\hat{x}_n$ (since $\hat{x}_n>\alpha_n$, the unique fixed point in $(0,1)$). Note also that in the case
 of existence of $a_{n}$ and $b_{n}$, we have $0\leq a_n<\alpha_n<\hat x_n<b_n<1$. We have moreover in this case the following
\begin{Lem}\label{minmax}
Assume the existence of $a_n$ and $b_n$  and  that $k:=\max\big(\left| f_{n}^\prime(a_n)|, |f_{n}^\prime(b_n)\right|\big)<1$. Then, the basin of  attraction of $\alpha_n$ is $[0,1)$.  
\end{Lem}

\begin{proof}
For typographical simplicity we chose to not write the subscripts $_{n}$.\\
Note that $f([a,b])=[f(\hat x),\hat x]\subset [f(\hat x),b]$ since $f$ is decreasing
on $[a,\hat x]$ and increasing on $[\hat x,b]$.
Moreover, using the convexity of $f$, first:
$$
\hat x - f(\hat x) = f(a) - f(\hat x) \leq \left|f^\prime(a)\right|(\hat x-a)\leq \hat x-a\,,
$$
implying that $f(\hat x) \geq a$ and thus $f([a,b])\subset[a,b]$,
and secondly: for every
$x\in[a,b]$,
\[\left|f(x)-\alpha\right|=\left|f(x)-f(\alpha)\right|\le \max_{a\leq y\leq b}\left|f^\prime(y)\right||x-\alpha| \leq k |x-\alpha|, \]
implying 
that the basin of attraction of $\alpha$ contains $[a,b]$.\\
Now, if $x\in(b,1)$, there exists $N\in\mathbb N^*$ such that $f^N(x)\le b$ and for all $0\leq q<N,\, f^q(x)\in(b,1)$. Indeed, since
  $f(y)<y$ on $(b,1)$, if $N$ does not exist, the sequence $(f^q(x))_{q\geq 0}$ is decreasing and bounded below by $b$, so tends toward a fixed point of $f$ in $[b,1)$, which raises a contradiction.\\
 As a result, the definition of $N$ and the monotonicity of $f$ on $(b,1)$  imply
 \[a\le \hat x=f(b)\le f^N(x)\le b\]
 and  we conclude that $(b,1)$ is included in the  basin of attraction of $\alpha$ since  $[a,b]$ is.\\
 In particular, the basin of attraction of $\alpha$  contains $[a,1)$ and, as $f$ is decreasing on $[0,a)$, $f([0,a))=\left(\hat x, f(0)\right]\subset [a, 1)$, implying that
 $[0,a)$ and thus $[0,1)$ is included in the basin of attraction of $\alpha$.
\end{proof}

The rest of the proof of Theorem \ref{only} is obtained using computer assistance to find good approximations for the quantities $\hat x_n$, $f_n(\hat x_n)$, $\alpha_n$, $f_n^\prime(\alpha_n)$, $a_n$, $b_n$, $f_n^\prime(a_n)$ and $f_n^\prime(b_n)$: 
 for every even $4\leq n\le 26$, $a_n$ and $b_n$ exist and satisfy the assumptions of Lemma \ref{minmax}, see the following Table. Proposition~\ref{nonunif_CV}  thus implies the statement of Theorem~\ref{only} in this case.\\[0.2cm]

\begin{filecontents*}{tableau9.csv}
n	tmin	fn(tmin)	tfixe	fn_prime(tfixe)	t1	t2	fn_prime(t1)	fn_prime(t2)
4	0,2531	0,228	0,2288	-0,0659	0,1308	0,4264	-0,4724	0,2431
6	0,207	0,1818	0,1825	-0,0674	0,0936	0,414	-0,5641	0,1894
8	0,1766	0,1548	0,1554	-0,0595	0,0773	0,3889	-0,5877	0,1519
10	0,1547	0,1369	0,1373	-0,05	0,0685	0,359	-0,576	0,1254
12	0,1382	0,1238	0,1241	-0,0408	0,0631	0,3273	-0,5429	0,1059
14	0,1252	0,1138	0,114	-0,0324	0,0597	0,2948	-0,4959	0,091
16	0,1146	0,1059	0,106	-0,0251	0,0576	0,2621	-0,4392	0,0793
18	0,1059	0,0994	0,0994	-0,0187	0,0566	0,2296	-0,3753	0,0698
20	0,0985	0,0939	0,0939	-0,0131	0,0565	0,1973	-0,3055	0,0621
22	0,0922	0,0892	0,0892	-0,0083	0,0577	0,1652	-0,2301	0,0552
24	0,0867	0,0852	0,0852	-0,0042	0,0607	0,133	-0,1474	0,0476
26	0,0818	0,0816	0,0816	-0,0007	0,0702	0,0967	-0,0452	0,0278
28	0,0776	0,0785	0,0785	0,0024				
30	0,0738	0,0757	0,0757	0,0051				
\end{filecontents*}
\csvreader[separator=tab,
    tabular=|c|c|c|c|c|c|c|c|c|,
    table head=\hline\rowcolor{gray!99}\color{white}$n$ & \color{white}$\hat x_n$ & \color{white}$f_n(\hat x_n)$ & \color{white}$\alpha_n$ & \color{white}$f_n^\prime(\alpha_n)$ & \color{white} $a_n$ &\color{white} $b_n$ &\color{white} $f_n^\prime(a_n)$ &\color{white} $f_n^\prime(b_n)$,
    late after head=\\\hline,
    late after line=\\\hline]
    {tableau9.csv}{ n=\cola, tmin=\colb, fn(tmin)=\colc, tfixe=\cold, fn_prime(tfixe)=\cole, t1=\colf, t2=\colg, fn_prime(t1)=\colh, fn_prime(t2)=\coli}
    {\cola & \colb & \colc & \cold & \cole & \colf & \colg & \colh & \coli }
    
    \begin{Rem}
  \begin{enumerate} 
\item Note that for $n\in\lbrace28,30\rbrace$, the cells corresponding to $a_n$ and $b_n$ are empty since no pre-image of $\hat x_n$ exists in these cases. 
\item The strategy in this section can not be used for the case of a GW, whereas Lemma \ref{le.Budan} implies the non-repulsivity of the fixed point which is a big step to achieve our goal if we succeed to prove the unicity of the fixed point.  
\end{enumerate}
\end{Rem}

\subsection{Estimates on the fixed points \boldmath$\alpha_{n}$}
In this section we obtain bounds for the fixed points of $f_n$ depending on $n$. As previously, we denote for $n\in\N$,  $\xi_{2n}=2^{-2n}\binom{2n}{n}$.
\begin{Pro}\label{Pierre'sfavorite}
We have:
\begin{equation}
\forall\,n\geq 536\,,\ \   \frac{1}{\sqrt{2\pi(n+\frac14)}}<\xi_{4n}\le  \alpha_n
 \quad\text{and}\quad 
\forall\,n\geq 3\,,\ \  \alpha_n\le \xi_{n^{\#}}
\le \frac{\sqrt 2}{\sqrt{\pi (n-1)}}\,,
\label{estim}
\end{equation}
where $n^{\#}:=2\left[\frac{n}{2}\right]$. 
\end{Pro}

\begin{proof}
According to (see \eqref{Nass1})
\begin{equation}
\label{eq.Nass}
\forall n \geq 1\,,\ \ \ 
{   \frac{2}{\sqrt{2\pi(2n+1)}}}< \xi_{2n}<
{ \frac{1}{\sqrt{\pi n}}},
\end{equation}
we have just to prove lower bound
$\xi_{4n}\le  \alpha_n$
for $n$ large enough and the upper bound 
$\alpha_n\le \xi_{n^{\#}}$ for~$n\geq 3$.

Using moreover the monotonicity of the sequence $\left(\xi_{2n}\right)$ and 
\eqref{eq.binom-easy},
we have for every $n\geq 2$:
%

\begin{eqnarray*}
\displaystyle
f_n(\xi_{4n})&=&\sum_{k,0\le 2k\le n}\binom{n}{2k}\xi_{2k}\left(1-\xi_{4n}\right)^{2k}\xi_{4n}^{n-2k}\\
&\ge&\sum_{k,0\le 2k\le n}\binom{n}{2k}\xi_{n^{\#}-2}\left(1-\xi_{4n}\right)^{2k}\xi_{4n}^{n-2k}+(\xi_{n^{\#}}-\xi_{n^{\#}-2})\binom{n}{n^{\#}}\left(1-\xi_{4n}\right)^{n^{\#}}\xi_{4n}^{n-n^{\#}}\\
&\ge&\xi_{n^{\#}-2}\frac{1+(-1)^n(1-2\xi_{4n})^n-2\frac{n}{n-1}\xi_{4n}\left(1-\xi_{4n}\right)^{n^{\#}}}{2}\\
&=& \xi_{{n}^{\#}}\frac{{n^{\#}}}{{n^{\#}}-1}\frac{1+(-1)^n(1-2\xi_{4n})^n-2\frac{n}{n-1}\xi_{4n}\left(1-\xi_{4n}\right)^{n^{\#}}}{2}\\
&\ge&  \xi_{4n}\frac{{n^{\#}}}{{n^{\#}}-1}\frac{1+(-1)^n(1-2\xi_{4n})^n-2\frac{n}{n-1}\xi_{4n}\left(1-\xi_{4n}\right)^{n^{\#}}}{e^{\frac{3}{4n}}} \,,
\end{eqnarray*}	
where the last inequality arises from using twice  (see \eqref{Nass2})
\begin{equation}
\label{eq.Nass2}
\forall n \geq 1\,,\ \ \ 
\xi_{4n}<\frac{1}{\sqrt{2}}e^{\frac{1}{4n}}\xi_{2n}\,.
\end{equation}

Using now the relations
$e^{x}<1+\frac43 x$ for $x\le \frac12$, $(1-x)^n\le e^{-nx}$ for $x\in[0,1]$, and
 \eqref{eq.Nass}, which implies
${ \frac{1}{\sqrt{2\pi (n+\frac14)}}<\xi_{4n}<\frac{1}{\sqrt{2\pi n}} }$, we have:
\begin{align*}
\frac{1+(-1)^n(1-2\xi_{4n})^n-2{ \frac{n}{n-1}}\xi_{4n}\left(1-\xi_{4n}\right)^{n^{\#}}}{e^{\frac{3}{4n}}}&\ge\frac{1-(1-2\xi_{4n})^n-\frac{2}{\sqrt{2\pi n}}{ \frac{n}{n-1}}\left(1-\xi_{4n}\right)^{n^{\#}}}{\frac{ n+1}{ n }}\\
&\ge \frac{1-e^{-\frac{2n}{\sqrt{2\pi (n+\frac14)}}}-\frac{2}{\sqrt{2\pi n}}{ \frac{n}{n-1}}
e^{-\frac{n-1}{\sqrt{2\pi (n+\frac14)}}}}{\frac{ n+1}{ n }}
\end{align*}
Finally we can  check that 
$$
\frac{{n^{\#}}}{{n^{\#}}-1}\left(1-e^{-\frac{2n}{\sqrt{2\pi (n+\frac14)}}}-\frac{2}{\sqrt{2\pi n}}{ \frac{n}{n-1}}
e^{-\frac{n-1}{\sqrt{2\pi (n+\frac14)}}}
\right)>{\frac{n+1}{n }}
$$
for all $n$ sufficiently large, and computer assisted calculations show that ${ n\ge 536}$ is sufficient.

Hence,
we have
$f_n(\xi_{4n})\ge \xi_{4n}$ and thus $\alpha_n\ge \xi_{4n}$
for every $n\geq 536$.
\medskip

To obtain the upper bound of \eqref{estim}, let us write for $n\geq 2$:
\[f_n(\xi_{n^{\#}})=\sum_{k,0\le 2k\le n}\binom{n}{2k}\xi_{2k}\left(1-\xi_{{n^{\#}}}\right)^{2k}\xi_{{n^{\#}}}^{n-2k}=:\sum_{k,0\le2k\le n}\nu_k\alpha_k\,,\]
where $\alpha_k:=\xi_{2k}$, and let $\mu_k:=\binom{n}{2k}$. The positive  sequence $(\alpha_k)_{k\geq 0}$ is decreasing according to Lemma \ref{XI}, and writing 
\[\frac{\nu_k}{\mu_k}=\left(\frac{1}{\xi_{{n^{\#}}}}-1\right)^{2k}\xi_{{n^{\#}}}^{n}, \]
the positive sequence $(\nicefrac{\nu_k}{\mu_k})_{k\geq 0}$ is increasing. Then, Lemma~\ref{wiard1} and the 
formulas \eqref{de.xi-n'}, \eqref{eq.binom-easy} give:  
%
\begin{align*}
f_n(\xi_{n^{\#}})
\le \frac{\sum_{k,0\le 2k\le n}\mu_k\alpha_k}{\sum_{k,0\le 2k\le n}\mu_k}\sum_{k,0\le 2k\le n}\nu_k
=\frac{2^{-n}\binom{2n}{n}}{2^{n-1}}\frac{1+(-1)^n(1-2\xi_{n^{\#}})^n}{2}\le \xi_{2n}\left(1+(1-2\xi_{n^{\#}})^n\right).
\end{align*}
Consequently, using in addition and $(1-x)^n\le e^{-nx}$ for all $x\in[0,1]$ and  the lower bound in \eqref{eq.Nass}:
\[\left(1-2\xi_{n^{\#}}\right)^n\le e^{-2n\xi_{n^\#}}
\le
{ e^{-\frac{4n}{\sqrt{2\pi(n^{\#}+1)}}}  }.
\]
Then, with \eqref{eq.Nass2}:
\[f_n(\xi_{n^{\#}})\le\xi_{2n^{\#}}\left(1+{  e^{-\frac{4n}{\sqrt{2\pi(n^{\#}+1)}}}  }\right) \le \xi_{n^\#}\frac{1}{\sqrt{2}} { 
e^{\frac{1}{2n^\#}}}\left(1+{  e^{-\frac{4n}{\sqrt{2\pi(n^{\#}+1)}}}  }\right)=: \xi_{n^\#}w_n .\]
Since $n^{\#}=n$ when $n$ is even and
$n^{\#}=n-1$ when $n$ is odd,
the sequences $(w_{2k})_{k\geq 1}$
and $(w_{2k+1})_{k\geq 1}$  are clearly decreasing and, as ${  w_3,w_{4}\le 1}$, we have
for every $n\geq 3$ :  $f_n(\xi_{n^{\#}})\le \xi_{n^{\#}}$ and thus $\alpha_n\le \xi_{n^\#}$.
\end{proof}

%

\section{An Example of GW}
\label{se.4}
All the simulations with a GW seem to show that there is a unique fixed point in $(0,1)$ and its basin of attraction is $(0,1)$. As we have already said, we have not been able to adapt the techniques of Section~\ref{se.3} to prove the uniqueness of the fixed point in a general framework. Nevertheless, we propose an example in which we are able to prove everything.\\
In this section, we  assume that the reproduction law $N$ follows a shifted geometric distribution with parameter $p\in(0,1)$, in other words:
\[q_n=\p(N=n)=p(1-p)^{n-2}, \forall n\ge2.\]
 This example is very satisfying as we can obtain explicit formulas.
More precisely, we have the following:

\begin{Lem}
If $N=X+1$ where $X$ follows a geometric distribution with parameter $p\in(0,1)$, we have: 
\begin{equation}
f(t)=\frac{p}{(1-p)^2}\left(-((1-p)t+1)+\frac{1}{(p(2-p+2t(p-1)))^{\frac{1}{2}}}\right).
\end{equation}
\end{Lem}

\begin{figure}[h!]
\begin{center}
\includegraphics[width=7cm]{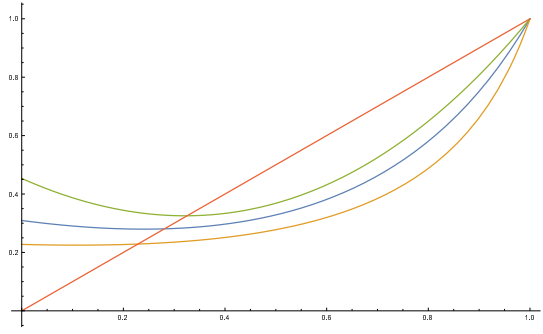}
\end{center}
\caption{The graphs of $f$ for $p=\nicefrac{1}{2}$ (blue), $p=\nicefrac{1}{4}$ (orange), $p=\nicefrac{9}{10}$ (green) }
\label{geom}
\end{figure}

\begin{proof}
We have:
\begin{align*}
f(t)&=\sum_{n\ge2}q_nf_n(t)=\sum_{n\geq 2} p(1-p)^{n-2}\frac{1}{\pi}\int_{0}^\pi((1-t)\cos x+t)^n\mathrm{d}x\\
&=\frac{p}{\pi}\int_{0}^\pi((1-t)\cos x+t)^2\sum_{n\geq 2} (1-p)^{n-2}((1-t)\cos x+t)^{n-2}\mathrm{d}x\\
&=\frac{p}{\pi}\int_{0}^\pi\frac{((1-t)\cos x+t)^2}{1- (1-p)((1-t)\cos x+t)}\mathrm{d}x=\frac{p}{\pi(1-p)^2}\int_{0}^\pi\frac{(1-p)^2((1-t)\cos x+t)^2-1+1}{1- (1-p)((1-t)\cos x+t)}\mathrm{d}x\\
&=\frac{p}{\pi(1-p)^2}\int_{0}^\pi-(1-p)((1-t)\cos x+t)-1+\frac{1}{1- (1-p)((1-t)\cos x+t)}\mathrm{d}x\\
&=\frac{-p((1-p)t+1)}{(1-p)^2}+\frac{p}{\pi(1-p)^2}\int_{0}^\pi\frac{\mathrm{d}x}{1- (1-p)((1-t)\cos x+t)}
\end{align*}
With the substitution $u=\tan\frac{x}{2}$, we obtain:
\begin{align*}
\int_{0}^\pi\frac{\mathrm{d}x}{1- (1-p)((1-t)\cos x+t)}&=2\int_0^{+\infty}\frac{\mathrm{d}u}{p+(2-p+2t(p-1))u^2}=\frac{2}{p}\int_0^{+\infty}\frac{\mathrm{d}u}{1+\frac{(2-p+2t(p-1))}{p}u^2}\\
&=\frac{\pi}{(p(2-p+2t(p-1)))^{\frac{1}{2}}}.
\end{align*}
Then
\[f(t)=\frac{p}{(1-p)^2}\left(-((1-p)t+1)+\frac{1}{(p(2-p+2t(p-1)))^{\frac{1}{2}}}\right).\]
\end{proof}

In Figure \ref{geom}, we can see that $f$ seems to have one fixed point on $(0,1)$. It is not difficult to find this point, resolving:  
\begin{align*}
f(t)=t&\Leftrightarrow -2(1-p)t^3+(2-5p)t^2+4pt-p=0\Leftrightarrow (t-1)(-2(1-p)t^2-3pt+p)=0\\
&\Leftrightarrow t=1\mbox{ or } t=\frac{-3p\pm(p(p+8))^{\frac{1}{2}}}{4(1-p)}.
\end{align*} 
And we can easily  see that  the only root that interests us is $\alpha=\frac{-3p+(p(p+8))^{\frac{1}{2}}}{4(1-p)}.$

\begin{Lem}
The basin of attraction  of $\alpha$ is $[0,1)$. 
\end{Lem}

\begin{proof}
The first and second derivative of $f$ are given by:
\[f^\prime(t)=\frac{p}{1-p}\left(-1+\frac{p}{(p(2-p+2t(p-1))^{\frac{3}{2}}}\right) \mbox{ and } f^{(2)}(t)=\frac{3p^3}{(p(2-p+2t(p-1))^{\frac{3}{2}}}.\]
As stated in Remark \ref{remark.1} since $f$ is strictly convex,   it is sufficient to show that $f^\prime(\alpha)\ge0$. 
One can see that:
\begin{align}\label{geomattract}
f^\prime(\alpha)\ge 0&\Leftrightarrow \frac{p}{1-p}\left(-1+\frac{p}{(p(2-p-\frac{1}{2}(-3p+(p(p+8))^{\frac{1}{2}})))^{\frac{3}{2}}}\right)\ge 0\nonumber\\
&\Leftrightarrow 2^{\frac{3}{2}}p\ge (p(4+p-(p(p+8))^{\frac{1}{2}})))^{\frac{3}{2}}\nonumber\\
&\Leftrightarrow 2\ge p^{\frac{1}{3}} (4+p-(p(p+8))^{\frac{1}{2}}=:g(p).
\end{align}
As $g(1)=2$, if we prove that $g$ is an increasing function on $[0,1]$, we obtain formula \eqref{geomattract}. As:
\begin{align*}
g^\prime(p)=\frac{(4+p-(p(p+8))^{\frac{1}{2}})((p(p+8))^{\frac{1}{2}}-3p)}{3p^{\frac{2}{3}}(p(p+8))^{\frac{1}{2}}}
\end{align*}
is obviously positive, our proof is complete.
\end{proof}

To conclude, according to Proposition \ref{nonunif_CV}, for every $k\geq 2$ and $\bp\in \mathscr P_k$ such that $\bp_1=\bp_2>\bp_{3}\geq\dots\geq \bp_k$, $\bp(m)$ converges to $\left( \alpha,\frac{1-\alpha}{2},\frac{1- \alpha}{2},0_{k-2}\right)$ when $m\to\infty$.

\begin{Rem}
Considering the $n$-ary tree as a GW tree with reproduction law $N=n$ a.s., we have~$\E[N]=n$. In order to obtain the same mean in the geometric case, it suffices to take~$p=\frac{1}{n-1}$.  With this choice of $p$, $\alpha\sim \nicefrac{1}{\sqrt{2n}}$ when $n$ goes to infinity, which is consistent with the bounds found in Proposition \ref{Pierre'sfavorite}.
\end{Rem}


\section{Open questions and variant case}
\label{se.5}

As a conclusion we make some remarks on the properties of the main objects studied in this work and discuss  possible generalisations of our results.
\begin{enumerate}
\item One can notice in the figure \ref{infcut}, that the red curve of $f_3$, seems to cut the blue one $f_4$ at its minimum. In fact, that is true for all $n\ge2$, that is:
\begin{equation}
f_{2n}\left(\hat x_{2n}\right)=f_{2n-1}\left(\hat x_{2n}\right).\label{argmin}
\end{equation}
Indeed, according to \eqref{deriv}:
\begin{align}
f_n^\prime(t)
&=\frac{n}{\pi}\int_0^\pi(t(1-\cos x)+\cos x)^{n-1}\mathrm{d}x-\frac{n}{\pi}\int_0^\pi\cos x(t(1-\cos x)+\cos x)^{n-1}\mathrm{d}x\nonumber\\
&=nf^\prime_{n-1}(t)-\frac{n}{\pi}\int_0^\pi\cos x(\cos x(1-t)+t)^{n-1}\mathrm{d}x\nonumber\\
&=nf^\prime_{n-1}(t)-\frac{n}{\pi}\int_0^\pi\frac{\cos x(1-t)+t-t}{1-t}(\cos x(1-t)+t)^{n-1}\mathrm{d}x\nonumber\\
&=n\left(1+\frac{t}{1-t}\right)f_{n-1}(t)-\frac{n}{1-t}f_n(t)\Leftrightarrow \frac{t-1}{n}f_n^\prime(t)=f_n(t)-f_{n-1}(t)\label{strange}
\end{align}
and taking $t=\hat x_{2n}$:
\[\frac{\hat x_{2n}-1}{n}f_{2n}^\prime(\hat x_{2n})=0=f_{2n}(\hat x_{2n})-f_{2n-1}(\hat x_{2n}).\]
With an obvious induction reasoning, the formula \eqref{strange} gives: 
\begin{equation}
f_{n}(t)-t=f_n(t)-f_{1}(t)=(t-1)\sum_{k={ 2}}^n\frac{1}{k}f_k^\prime (t).\label{telesco}
\end{equation}

We think that these equalities have a probabilistic meaning, but did not manage to come up with an explanation.

\item The general GW case for two opinions seems for the moment out of reach, even though our  simulations suggest that  our results stay valid.
Contrary to the  article \cite{Haber}, the mean of the reproduction law $\E[N]$ does not seem to play a particular role: there seems to be always convergence. 
\\
 We have to study cases with more than two opinions: nevertheless, even in the case of a $n$-ary tree, using links with random walks in order to obtain formulas like \eqref{Luc} for  a number $i>2$ of major opinions is not clear to us. Moreover, we have seen that even in the case with two opinions,  parity plays an important role; already for $i=3$, the calculations become devilish and it seems to us that we need to find much finer methods than direct computations.\\
For instance in Figure \ref{third}, we can see that the shape of the graph is linked to the remainder of the Euclidean division of $n$ by 3 and the equivalent formula for $f_n$ is: 
\begin{align*}
f_n(t)&=\sum_{k=0}^n\binom{n}{3k}\binom{3k}{k}\binom{2k}{k}t^{n-3k}\left(\frac{1-t}{3}\right)^{3k}\\
&+3\sum_{k=0}^n\binom{n}{2k}\binom{2k}{k}\sum_{j=0}^{n-2k}\binom{n-2k}{j}t^{n-2k-j}\left(\frac{1-t}{3}\right)^{2k+j}.
\end{align*}
\begin{figure}[h!]
\begin{center}
\includegraphics[width=7cm]{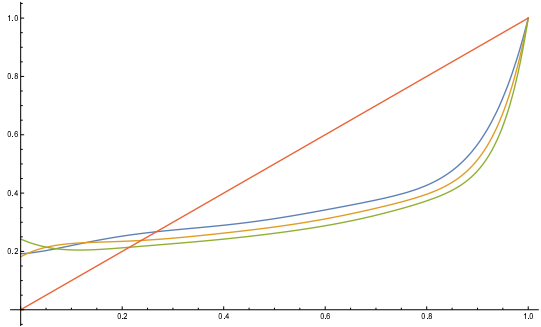}
\end{center}
\caption{The graphs of $f_n$ for $n=7$ (blue), $n=9$ (orange), $n=11$ (green) }
\label{third}
\end{figure}

\end{enumerate}

\section{Appendix}
\label{se.6}

In this appendix, we recall some classical definitions and results used throughout this paper.\\
The following result is crucial to prove the stability statement of Section~\ref{sub.even-attract}.


\begin{Lem}\label{wiard1}
Consider two positive sequences $(\mu_k)_{k\geq 0}$ and $(\nu_k)_{k\geq0}$ such that $\left(\nicefrac{\nu_k}{\mu_k}\right)_{k\geq 0}$ is increasing.  
Then, for every decreasing (resp. increasing) sequence $(\alpha_k)_{k\geq 0}$ and for every 
$0\leq \ell\leq n$,
\begin{equation}\label{metalemma}
\frac{\sum_{k=\ell}^{n}\alpha_k\mu_k}{\sum_{k=\ell}^{n}\mu_k}\ \ge\  \frac{\sum_{k=\ell}^{n}\alpha_k\nu_k}{\sum_{k=\ell}^{n}\nu_k}
\ \ \Big(\text{resp.}\ \le\  \frac{\sum_{k=\ell}^{n}\alpha_k\nu_k}{\sum_{k=\ell}^{n}\nu_k}\ \Big)\,.
\end{equation}
\end{Lem}

\begin{proof}
Assume that the sequence $\left(\nicefrac{\nu_k}{\mu_k}\right)_{k\geq 0}$ is increasing, which is equivalent to: 
\[\forall\,0\le i\le j\,,\ \ \  {\mu_i\nu_j}\ge {\mu_j\nu_i}.\]
When the sequence $(\alpha_k)_{k\geq 0}$ is decreasing, the  formula \eqref{metalemma} 
 is equivalent to: 
\begin{align*}
 \sum_{i,j=\ell}^{n}\mu_i\alpha_i\nu_j\ge \sum_{i,j=\ell}^{n}\mu_i\alpha_j\nu_j
&\Leftrightarrow \sum_{\ell \le i<j \le n}\mu_i\alpha_i\nu_j+\sum_{\ell \le j\le i \le n}\mu_i\alpha_i\nu_j\ge
\sum_{\ell \le i<j \le n}\mu_i\alpha_j\nu_j+\sum_{\ell \le j\le i \le n}\mu_i\alpha_j\nu_j\\
&\Leftrightarrow \sum_{\ell \le i<j \le n}\mu_i\alpha_i\nu_j+\sum_{\ell \le i\le j \le n}\mu_j\alpha_j\nu_i\ge\sum_{\ell \le i<j \le n}\mu_i\alpha_j\nu_j+\sum_{\ell \le i\le j \le n}\mu_j\alpha_i\nu_i\\
&\Leftrightarrow \sum_{\ell \le i<j \le n}\mu_i\alpha_i\nu_j+\sum_{\ell \le i<j \le n}\mu_j\alpha_j\nu_i\ge\sum_{\ell \le i<j \le n}\mu_i\alpha_j\nu_j+\sum_{\ell \le i<j \le n}\mu_j\alpha_i\nu_i\\
&\Leftrightarrow \sum_{\ell \le i<j \le n}(\mu_i\nu_j-\mu_j\nu_i)(\alpha_i-\alpha_j)\ge 0\,,
\end{align*}
which is true by hypothesis.
\end{proof}

\ 


In the following lemma, we state
classical results on binomial coefficients.

\begin{Lem}\label{wiard3}
For all $n\in\mathbb N$: 
\begin{align}
\sum_{k=0}^n (2k+1)\binom{n}{2k+1}=\sum_{k=0}^n 2k\binom{n}{2k}&=n2^{n-2}\ \ \text{ when $n\neq 1$}\label{combineasyimpli}\\
\text{and}\qquad  \sum_{j=0}^{n}\binom{n}{2j}\binom{2j}{j}2^{-2j}&=2^{-n}\binom{2n}{n}\,.\label{luc1}
 \end{align}
\end{Lem}
\begin{proof}
\begin{enumerate}
\item For all $x,y\in\mathbb R$ and $n\geq 0$, the relation
\begin{equation}
\label{eq.binom-easy}
h(x,y):=\sum_{k=0}^n \binom{n}{2k}x^{2k}y^{n-2k}=\frac{(x+y)^n+(y-x)^n}{2}
\end{equation}
implies that
\[\frac{\partial h }{\partial x}(x,y)=\sum_{k=0}^n 2k\binom{n}{2k}x^{2k-1}y^{n-2k}=\frac{n(x+y)^{n-1}-n(y-x)^{n-1}}{2}.\]
Taking $x=y=1$ and $n\neq 1$, we obtain the right equality of \eqref{combineasyimpli}.\\
Moreover, with a very similar reasoning:
\[g(x,y)=(x+y)^n=\sum_{k=0}^n\binom{n}{k}x^ky^{n-k}\mbox{ and }\frac{\partial g}{\partial x}(x,y)=n(x+y)^{n-1}=\sum_{k=0}^n\binom{n}{k}kx^{k-1}y^{n-k}\]
and taking $x=y=1$, we obtain 
\[\sum_{k=0}^n\binom{n}{k}k=n2^{n-1}.\]
We conclude by using
\[\sum_{k=0}^nk\binom{n}{k}=\sum_{k=0}^n 2k\binom{n}{2k}+\sum_{k=0}^n (2k+1)\binom{n}{2k+1}.\]
\item In \cite{Gould}, the author uses an expansion of $(x^2+2x)^n$ to prove \eqref{luc1} (see the formula  (1.65) there). 
We will use here the following series expansion,
 for $\ell\in\mathbb N$ and $x\in[0,1)$ (which also permit to prove the generalization of \eqref{luc1}  stated in Remark~\ref{re.gen} below):
\begin{equation*}
\frac{x^\ell}{(1-x)^{\ell+1}}=\sum_{n\ge0}\binom{n}{\ell}x^n\quad\text{and}\quad
\frac{1}{\sqrt{1-x}}=\sum_{n\ge 0}\binom{2n}{n}2^{-2n}x^{n}.
\end{equation*}
The first one can be obtained by induction and the second one is classical. Thus:
\begin{eqnarray*}
\sum_{n\ge 0}\sum_{j\ge0}\binom{n}{2j}\binom{2j}{j}2^{-2j}x^n&=&\sum_{j\ge0}\binom{2j}{j}2^{-2j}\sum_{n\ge 0}\binom{n}{2j}x^n=\sum_{j\ge0}\binom{2j}{j}2^{-2j}\frac{x^{2j}}{(1-x)^{2j+1}}\\
&=&\frac{1}{1-x}\sum_{j\ge0}\binom{2j}{j}2^{-2j}\left(\frac{x}{1-x}\right)^{2j}=\frac{1}{1-x}\times \frac{1}{\sqrt{1-\left(\frac{x}{1-x}\right)^2}}\\
&=&\frac{1}{\sqrt{1-2x}}=\sum_{n\ge 0}\binom{2n}{n}2^{-2n}(2x)^{n}=\sum_{n\ge 0}\binom{2n}{n}2^{-n}x^{n}.
\end{eqnarray*}
Identifying the coefficients, we obtain \eqref{luc1}.
\end{enumerate}
\end{proof}

\begin{Rem}
\label{re.gen}
Adapting the above proof of \eqref{luc1}, we can show that for all $n\in\mathbb N$ and all $\ell \in\mathbb N^*$:
\begin{equation*}
\sum_{j=0}^{n}j\dots(j-\ell+1)\binom{n}{2j}\binom{2j}{j}2^{-2j}=2^{-n}\binom{2(n-\ell)}{n-\ell}(n-\ell)\dots (n-2\ell+1)\label{genfor}.
\end{equation*}
\end{Rem}

%


We conclude this appendix with
this last  lemma,  giving some properties of the Wallis integrals and of the strongly related quantities $\xi_{2n}=\frac{1}{2^{2n}}\binom{2n}{n}$
defined in \eqref{de.xi-n}.

\begin{Lem}
\label{XI}
For all $n\in\mathbb N$, define the  Wallis integral
\begin{equation*}
W_n:=\int_0^{\nicefrac{\pi}{2}}\sin^nt\,\mathrm{d}t.
\end{equation*}
The sequence $(W_n)_{n\ge0}$ is positive and strictly decreasing and, for all $n\ge 1$:
\begin{equation}\label{easystuff}
\sqrt{\frac{\pi}{2(n+1)}} < W_n < \sqrt{\frac{\pi}{2n}}
\end{equation}  
We have moreover the following properties:
\begin{enumerate}
\item The sequence $(\xi_{2n})_{n\ge0}$ is strictly decreasing.
\item For all $n\ge1$: 
\begin{equation}
{   \frac{2}{\sqrt{2\pi(2n+1)}}}< \xi_{2n}<
{ \frac{1}{\sqrt{\pi n}}}\label{Nass1}
\end{equation}
{ 
and hence
\begin{equation}
\label{Nass2}
\xi_{4n}<\frac{1}{\sqrt{2}}e^{\frac{1}{4n}}\xi_{2n}\,.
\end{equation}
}

\end{enumerate}
\end{Lem}

\begin{proof}
Let us prove the well-known formula \eqref{easystuff} for the sake of completeness.
For $n\ge0$, as $0\le \sin t\le 1$ in $[0,\nicefrac{\pi}{2}]$ (and $0< \sin t< 1$ in $(0,\nicefrac{\pi}{2})$): 
\begin{align*}
W_{n}>0\ \ \text{and}\ \ W_{n+1}-W_n=\int_{0}^{\nicefrac{\pi}{2}}\sin^n t(\sin t-1)\mathrm{d}t<0, 
\end{align*}
implying the (strict) monotonicity of $(W_n)_{n\ge0}$. Moreover, for every $n\in\mathbb N$:
\begin{align*}
 W_{n+2}&=\int_{0}^{\nicefrac{\pi}{2}}\sin t\sin^{n+1} t\mathrm{d}t=(n+1)\int_{0}^{\nicefrac{\pi}{2}}\cos^2t\sin^n t\mathrm{d}t=(n+1)(W_n-W_{n+2})\nonumber\\
& \Leftrightarrow (n+2)W_{n+2}=(n+1)W_n\\
&\Leftrightarrow (n+2)W_{n+2}W_{n+1}=(n+1)W_{n+1}W_n.\nonumber
\end{align*}
Consequently, the sequence $\big((n+1)W_{n+1}W_n\big)_{n\geq0}$ is constant and then:
\begin{equation*}\label{Wallis1}
\forall n\in\mathbb N^{*}\,,\ \ nW_{n}W_{n-1}=W_1W_0=\frac{\pi}{2}.
\end{equation*}
Using the monotonicity of $(W_n)$, we obtain the formula \eqref{easystuff} since:  
\begin{equation*}\label{Wallis2}
\forall n\geq 0\,,\ \ nW_n^2< \frac{\pi}{2}<(n+1)W^{2}_n.
\end{equation*}
Recall now that for all $n\ge 0$ (see \eqref{Luc1}),
\begin{equation*}
\label{evidence}
W_{2n}=\frac{\pi}{2}\frac{(2n)!}{2^{2n}(n!)^2}=\frac{\pi}{2}\xi_{2n}.
\end{equation*}
Thus:
\begin{enumerate}
\item The (strict) monotonicity of the sequence $(\xi_{2n})_{n\ge0}$ follows from  the one of $(W_n)_{n\ge0}$. 
\item Using \eqref{easystuff}, we obtain:
\[\sqrt{\frac{\pi}{2(2n+1)}}< W_{2n}< \sqrt{\frac{\pi}{4n}}\ \ \Longrightarrow\ \  \frac{2}{\sqrt{2\pi(2n+1)}}< \xi_{2n}< \frac{1}{\sqrt{\pi n }}.  \]
\end{enumerate}
\end{proof}

\section*{Acknowledgments}

We would like to thank our colleagues Luc Hillairet for his fruitful ideas and conversations, and Thomas Haberkorn for his help with simulations.
We also thank the anonymous referee who pointed out an argument permitting to state Corollary~\ref{Corcool}
for a GW tree with support included in $2\mathbb N+1$ without assuming its finiteness.

\end{document}